\documentclass[3p,letterpaper]{elsarticle}

\pdfoutput=1
\usepackage{times} 
\usepackage{amsmath,amssymb,amsfonts}
\usepackage[latin1]{inputenc}
 \usepackage[english]{babel}
\usepackage[cyr]{aeguill}
\usepackage{xcolor} 
\usepackage{geometry} 
\usepackage[pdftex,a4paper,hyperindex=true,hyperfigures=true,bookmarks=true,pdftitle={TITLE},pdfsubject={},pdfauthor={AUTHOR},pdfkeywords={KEYWORDS},pdfpagemode=UseOutlines,colorlinks=true,urlcolor=blue,citecolor=red,linkcolor=blue]{hyperref}
 \ifx\pdfoutput\undefined 
    \usepackage[dvips]{graphicx}
\else                    
 \usepackage[]{graphicx}
\fi
\usepackage{float} % Pour l'insertion d'images
\DeclareGraphicsExtensions{.jpg,.mps,.pdf,.png} % Formats d'images
\usepackage[sf]{caption2}
\usepackage[figuresright]{rotating}

\makeindex  % \index command in the .tex
\makeglossary  % \nomenclature command in .tex

\newcommand{\RR}{{\mathbb R}}

\newcommand{\ep}{{\varepsilon}}
\newcommand{\real}{{\mathbb R}}
\newcommand{\R}{{\mathbb R}}

\newcommand{\ZZ}{{\mathbb Z}}
\newcommand{\NN}{{\mathbb N}}

\newcommand{\la}{\lambda}

\newcommand{\al}{\alpha}
\newcommand{\be}{\beta}
\newcommand{\Hmin}{{h^{\min}}}
\newcommand{\pL}{\ell^{(p)}}
\newcommand{\pLt}{\ell^{(p,s)}} % fractional integral

\newcommand{\BE}{\begin{equation}}
\newcommand{\EE}{\end{equation}}

\newcommand{\ome}{\omega}

\newtheorem{lem}{Lemma}
\newtheorem{coro}{Corollary}
\newtheorem{Theo}{Theorem}
\newtheorem{prop}{Proposition}
\newtheorem{defi}{Definition}
\newcommand{\BP}{\begin{prop}}
\newcommand{\EP}{\end{prop}}
\newcommand{\BC}{\begin{coro}}
\newcommand{\EC}{\end{coro}}
\newcommand{\BL}{\begin{lem}}
\newcommand{\EL}{\end{lem}}
\newcommand{\BD}{\begin{defi}}
\newcommand{\ED}{\end{defi}}
\newcommand{\BT}{\begin{Theo}}
\newcommand{\ET}{\end{Theo}}

\newcommand{\review}[1]{{#1}}

\begin{document}
\begin{frontmatter}
\title{\bf\boldmath $p$-exponent and $p$-leaders, Part I:\\
Negative pointwise regularity.}
%\author{S. Jaffard, C. Melot, R. Leonarduzzi, H. Wendt,   P. Abry, S. G. Roux, M.~E. Torres} 
%\date{\today} 

\author[cret]{S. Jaffard}
\ead{jaffard@u-pec.fr}

\author[mars]{C. Melot}
\ead{melot@cmi.univ-mrs.fr}

\author[lyon]{R. Leonarduzzi}
\ead{roberto.leonarduzzi@ens-lyon.fr}

\author[toul]{H. Wendt}
\ead{herwig.wendt@irit.fr}

\author[lyon]{P. Abry\corref{mycorraut}}
\ead{patrice.abry@ens-lyon.fr}

\author[lyon]{S. G. Roux}
\ead{stephane.roux@ens-lyon.fr}

\author[arge]{M.~E. Torres}
\ead{metorres@santafe-conicet.gov.ar}

\address[cret]{%
Universit\'e Paris Est, Laboratoire d'Analyse et de Math\'ematiques Appliqu\'ees, CNRS UMR 8050, UPEC,  Cr\'eteil, France}
\address[mars]{%
Aix Marseille Université, CNRS, Centrale Marseille, I2M, UMR 7373, 13453 Marseille, France}  
\address[lyon]{%
Signal, Systems and Physics, Physics Dept., CNRS UMR 5672, Ecole Normale Sup\'erieure de Lyon, Lyon, France}
\address[toul]{%
IRIT, CNRS UMR 5505, University of Toulouse,  France}
\address[arge]{%
Consejo Nacional de Investigaciones Cient\'{\i}ficas y T\'ecnicas, Universidad Nacional de Entre Ríos, Argentina}

\cortext[mycorraut]{Corresponding author. Tel: +33 47272 8493. Postal address: CNRS, Laboratoire de Physique Ecole Normale Superieure de Lyon 46, allée d'Italie, F-69364, Lyon cedex 7, France}

\begin{abstract} 
Multifractal analysis aims to characterize signals, functions, images or fields, via the fluctuations  of their local regularity along time or space, hence capturing crucial features of their temporal/spatial dynamics.
Multifractal analysis is becoming a standard tool in signal and image processing, and is nowadays widely used in numerous applications of different natures. 
Its common formulation relies on the measure of local regularity via the H\"older exponent, by nature restricted to positive values, and thus to locally bounded functions or signals. 
It is here proposed to base the quantification of local regularity on $p$-exponents, a novel local regularity measure potentially taking negative values. 
First, the theoretical properties of $p$-exponents are studied in detail.
Second, wavelet-based multiscale quantities, the $p$-leaders, are constructed and shown to permit accurate practical estimation of $p$-exponents. 
Exploiting the potential dependence with $p$, it is also shown how the collection of $p$-exponents enriches the classification of locally singular behaviors in functions, signals or images.
The present contribution is complemented by a companion article developing the $p$-leader based multifractal formalism associated 
to $p$-exponents. 
\end{abstract}

%\keywords{keywords}
\end{frontmatter}
% \maketitle

\section{Introduction}

\noindent {\bf Context: multifractal analysis, local regularity and H\"older exponent.} \quad Multifractal analysis nowadays constitutes a classical signal and image processing tool, available in most modern toolboxes. 
It is commonly used to model or analyze scaling properties and local irregularities in signals and in image textures.
It has been involved in a large variety of real-world applications of very different natures, ranging from biomedical (heart rate variability \cite{Goldberger2002}, neurosciences \cite{Ciuciu2012,kantelhardt2002multifractal}), to physics (turbulence \cite{m74}), geophysics (rainfalls \cite{Foufoula-Georgiou1994}, \review{wind \cite{telesca2011analysis}, earthquakes \cite{telesca2006}}), finance \cite{Mandelbrot1997,Mandelbrot1999,Lux2007}, \review{music \cite{telesca2011revealing}}, or Internet traffic \cite{Abry2002}, to name but a few.

In essence, multifractal analysis describes, via the multifractal spectrum, the fluctuations of local regularity along time, or space, of a function, signal, image or field. 
Fundamentally, it relies on two key ingredients: 
At the theoretical level, a pointwise regularity exponent $h(t)$ that formalizes the intuition of local regularity;
at the practical level, a multiscale quantity $T_X(a,t)$ that permits the actual measurement of local regularity, via a local power-law behavior as a function of the analysis scales $a$: 
\begin{equation}
\label{loglog}
T_X(a,t)  \simeq a^{h(t)} \quad \mbox{when} \quad {a \rightarrow 0},
\end{equation}
or, more technically, 
\begin{equation}
\label{loglog}
h(t) = \liminf_{a \rightarrow 0} \frac{\log (T_X (a, k(t))) }{\log a }.
\end{equation}

In its commonly, if not exclusively, used formulation, multifractal analysis relies on the quantification of local regularity via the so-called H\"older exponent \cite{Fal93,jmf2,riedi03,Jaffard2004}.
In practice, several multiscale quantities have been involved in multifractal formalisms, the practical counterpart of multifractal analysis. 
It has now been long  recognized that increments or wavelet coefficients do not yield accurate analysis of H\"older exponents. 
The earliest solid practical formulation of multifractal analysis was based on the (continuous) wavelet transform modulus maxima method (WTMM) (cf. e.g., \cite{muzy1993multifractal,muzyetal94,ard02,Arneodo2003a}).
A multifractal extension of Detrended Fluctuation Analysis (MF-DFA) has also been proposed (cf. \cite{kantelhardt2002multifractal} for the founding article, \review{and \cite{gu2010detrending,schumann2011multifractal} for related developments}). 
Recently, a theoretically well-grounded and practically efficient formulation that, unlike WTMM and MFDFA, extends well to higher dimensional signals, has been proposed:
It relies on wavelet leaders, constructed as local suprema of discrete wavelet coefficients, cf. \cite{Jaffard2004,JAFFARD:2006:A,Wendt:2007:E,WENDT:2009:C}.\\

\noindent {\bf Limitation: negative regularity.} \quad An important, yet often overlooked, issue in multifractal analysis, consists in the fact that choosing multiscale quantities on which analysis is based  also amounts to selecting a specific measure of pointwise regularity. 
For instance, the wavelet leaders are the only multiscale quantities that were shown theoretically to be able to actually characterize the H\"older exponent \cite{Jaffard2004,JAFFARD:2006:A,Wendt:2007:E,WENDT:2009:C}.
Conversely, the practical use of a particular regularity exponent necessitates the construction of multiscale quantities specifically tailored to it, requiring the verification of global a priori regularity assumptions, which may not always hold for the data to be analyzed. 
In particular, by definition, H\"older exponents cannot take negative values.
In its current formulation, multifractal analysis is thus restricted to functions whose local regularity is everywhere positive, i.e., to functions that are everywhere locally bounded and thus show a minimal global regularity.
This is a severe restriction for the practical use of multifractal analysis as signals and images from real-world applications are often characterized by discontinuities and thus by a negative uniform, or minimum global, regularity. This is notably quasi-systematically the case for biomedical time series,
 (cf. e.g., \cite{Bergou,MandMemor} and references therein for reviews). 

It has often been proposed to circumvent this well-recognized limitation by integrating data prior to performing multifractal analysis \cite{muzy1993multifractal,muzyetal94,ard02,Arneodo2003a}).
The use of fractional integration, of order tuned to the uniform regularity of data, has recently been proposed and embedded directly in the definition of wavelet leaders \cite{WENDT:2009:C,Bergou,MandMemor}.
However, fractional integration yields the significant issue of relating the characterization of local regularity of integrated functions to those of the original functions, prior to integration. 
Notably, the commonly used rule of thumb that fractional integration results in a uniform shift of the multifractal spectrum is in general wrong (specifically, in the presence of oscillating singularities \cite{WENDT:2009:C}).
This issue will be further discussed in Section~\ref{sec:pExp}.
\\

\noindent {\bf Goals, contributions and outline.} \quad In this context, as a possible alternative to (fractional) integration, we propose the use of $p$-exponents, which potentially take negative values and hence permit to characterize negative local regularity.
Though introduced in the theoretical context of PDEs as early as 1961 for $p>1$ by Calder\'on and Zygmund \cite{CalZyg}, $p$-exponents were not used in signal processing until the 2000s when their wavelet characterization was proposed \cite{JaffMel,JaffCies,JaffToul}. 
In the present contribution, we study which information on the local behavior of a function near a singularity is supplied by the knowledge of the collection of $p$-exponents.

After recalling definitions, Section~\ref{sec:pwre} studies the theoretical properties of $p$-exponents. 

The corresponding multiscale quantities, the $p$-leaders, are defined in Section~\ref{sec:wavtool} as local $\ell^p$ norms of (discrete) wavelet transform coefficients. 
They will be shown to permit an accurate practical measure of $p$-exponents as well as to extend the definition of $p$-exponents to $p \in(0,1)$. 
Several representative examples will be used to illustrate the properties and potential of $p$-leaders and $p$-exponents to characterize local singular behaviors.

The joint use of all $p$-exponents paves the way to a thorough classification  of local singular behaviors (in terms of  $p$-invariant, canonical, oscillating, balanced or lacunary singularities); this will be detailed in Section \ref{sec:pExp}.
The interplay with negative regularity will also be discussed (Section~\ref{sec:negative}). 

While the present contribution addresses pointwise or local regularity and hence focuses on isolated singularities, the final goal is multifractal analysis, aiming to study collections of intertwined singularities.
Therefore, the present work is complemented by a companion article which develops and studies the corresponding $p$-leader based multifractal formalism that permits the practical implementation of multifractal analysis based on $p$-exponents. 

{\sc Matlab} routines designed by ourselves and implementing $p$-leaders and estimation procedures for $p$-exponents will be made publicly available to the research community at the time of publication of the present article.

\section{Pointwise  regularity exponents}
\label{sec:pwre}  

Throughout the present article, $\{ X(x), x \in \real^d \} $ denotes the function or sample path of a stochastic process, or random field, in dimension $d$, to analyze.  

\subsection{Pointwise H\"older regularity}
\label{sec:holex}

The most commonly used notion of pointwise regularity is defined via the so-called H\"older exponent, whose definition and properties are recalled below.
For Section~\ref{sec:holex}, $\{ X(x) \}_{x\in \R^d} $ is assumed to consist in a locally bounded function.\\

\noindent {\bf\boldmath Local H\"older spaces $C^{\alpha}(x_0)$.} \quad 
Let $\al \geq 0$.
The function $X$  is said to belong to $C^{\alpha}(x_0)$ at  location $x_0 \in \R^d$, with ${\alpha} \geq 0$, if there exist a constant $C >0$ and a polynomial $P_{x_0}$ of degree less than $ \alpha $, such that, for $a$ small enough,
\begin{equation}
\label{equ:tpe}
|X(x_0+a)-P_{x_0}(x_0+a)| \leq C |a|^{\alpha}. 
\end{equation}
Because it obviously generalizes the Taylor polynomial  for $C^N$ functions, $P_{x_0}$  is usually also referred to as the {\em Taylor polynomial}  of order $\al$  of $X$ at $x_0$. 
When  $\alpha < 1$, the Taylor polynomial boils down  to a constant $P_{x_0}(x) \equiv X(x_0)$. 
A more general discussion on $P_{x_0}$ will be given in Section \ref{sec:pexpodef} (see also \cite{Abel}).\\

\noindent {\bf Pointwise H\"older exponent.} \quad The \emph{H\"older exponent} of $X$ at location $x_0$ is defined as
\begin{equation}
\label{equ:he}
h(x_0)=\sup\{\alpha: X \in C^{\alpha}(x_0)\}.
\end{equation}
It is commonly used to characterize the local regularity of $X$ at $x_0$: The larger $h(x_0)$ the smoother $X$ around $x_0$. \\

\noindent {\bf Limitation of pointwise H\"older exponent: positive regularity.} \quad By definition, the H\"older exponent $h(x_0)$ cannot take negative values and (\ref{equ:tpe}) implies that $X$ is bounded in a neighborhood of $x_0$.
Indeed, $P_{x_0}$ is bounded and, since $\al \geq 0$,  so is $|a|^{\alpha}$.  
However, large classes of signals and images from real-world applications cannot be satisfactorily modeled by locally bounded functions (cf., e.g., \cite{MandMemor} where a practical  wavelet based criterium to assess local boundedness  is supplied and applied to numerous real-world data). 

The attempt to characterize singularities with negative regularity by simply allowing  $\al$ in  (\ref{equ:tpe}) to take negative values does not yield a satisfactory definition:  
Indeed, if (\ref{equ:tpe}) is satisfied for a given $\al$ (which may be negative), then, in any corona $0 < C_1 \leq | x-x_0 | \leq C_2$, $f$ is  a bounded function.
Therefore, such a definition would only permit to define {\em isolated singularities} of negative order, which would thus not be relevant in the framework of multifractal analysis where singularities with a given exponent may be dense, and usually consist of sets with strictly positive Hausdorff dimensions (whereas a set of isolated points is at most  countable and hence with Hausdorff dimensions equal to $0$).

\subsection{$p$-exponent regularity}
\label{sec:pexpodef}

To characterize negative regularity in data, a new  definition for  pointwise regularity is used, the $T^p_\alpha (x_0)$ regularity,  which has the advantage of relying on the less restrictive assumption that data  locally belong to  $ L^p (\R^d)$, 
instead of requiring local boundedness. 
This notion was introduced  by A. Calder\'on and A. Zygmund \cite{CalZyg} and has recently been put forward in the mathematical literature in \cite{JaffMel} (see also \cite{JaffToul,Bergou}). 
It permits the definition of a collection of $p$-exponents to measure pointwise regularity.

\BD \label{def:Tpal} Let $p \geq 1$ and $X$  be a function that locally belongs to $ L^p (\R^d)$.  
Let  $B(x_0, a)$ denote the ball centered at  $x_0$ and of radius $a >0$. 
Let  $\alpha > -d/p$.  
A function $ X $  belongs to $T^p_\alpha (x_0)$  if there exist a constant $ C$  and a polynomial $P_{x_0}$  of degree  less than $\alpha $  such that, for $a$  small enough,  
\begin{equation}  
\label{pexpa} 
T^{(p)}(a,x_0) : = \left( \frac{1}{a^d} \int_{B(0,a)} | X(u+x_0) -P_{x_0}(u+x_0)|^p du  \right)^{1/p} \leq C a^\alpha.
\end{equation}
%Note that, in 1D, $ B(x_0,a) $ is the interval $[x_0-a, 0+a]$.  \\
\ED
The \emph{$p$-exponent}  of  $X$ at $x_0$ is defined in \cite{JaffMel} by:
\begin{equation}  
\label{equ-pexp} 
h_p(x_0) = \sup \{ \alpha : X \in T^p_\alpha (x_0)\}.
\end{equation}
The $p$-exponent $h_p(x_0)$  takes values in $[ -d/p, + \infty ]$ (see Theorem \ref{th:tha} below), thus allowing for a proper mathematical definition of \emph{negative order pointwise regularity}. 
This definition constitutes a natural substitute for the pointwise H\"older regularity when dealing with functions which are not locally bounded.
The usual H\"older regularity actually corresponds to  $p=+\infty$ (in this case, the local $L^p$ norm in the left hand side of (\ref{pexpa}) boils down to a local $L^\infty$ norm, thus supplying a condition equivalent to (\ref{equ:tpe})). 
Definition \ref{def:Tpal} does not include $p$-exponents with $p <1$ because $L^p$ with $0<p<1$ is mathematically ill-defined.
Section \ref{sec:ppluspetun} will however  show how the wavelet framework permits to extend the definition of $p$-exponents to $0<p <1$. \\

\noindent {\bf Taylor polynomial.} \quad The polynomial $P_{x_0}$ defined implicitly by (\ref{pexpa}) is unique for a given $\al$ independently of $p$. 
When $\al$ crosses an integer value $N$, $P_{x_0}$  may be modified by addition of  terms of degree $N$. 
However, for two different $ \al_1$ and $\al_2$, the expansions of $P_{x_0}$ coincide up to order $\min (\al_1, \al_2)$ (cf. \ref{sec:prooftha}). 
This implies that by picking up the integer part of $h_p(x_0)$ we get the polynomial $P$ which corresponds to the largest possible value of $\al$. 
Theorem \ref{th:tha} (Point \ref{it:dec}) below indicates that one can fix a unique {\em Taylor polynomial} of $X$ at $x_0$, whose coefficients are independent of $p$, and are referred to as the (generalized) Peano derivatives  of $X$ at $x_0$ \cite{Ash,Abel}.

One of the main advantages of the wavelet framework developed in Section \ref{sec:wavtool} is however that the computation of $P_{x_0}$ is not  required to measure the $p$-exponent. 
The Taylor polynomial is thus not further discussed and readers are referred to e.g., \cite{Abel}.

\subsection{ Properties of $p$-exponents}
\label{sec:constucta}

Let us now state the main theorem characterizing the properties of the collection of $p$-exponents, or mapping $p \rightarrow h_p (x_0)$, at  location $x_0$.
\review{Let us also introduce the notation $X\in L^p (\RR^d) $ to indicate that $X$ belongs to $L^p (\RR^d)$ locally.}
 
\BT 
\label{th:tha}
Let $X\in L^1_{loc} (\RR^d) $ and  $x_0 \in \RR^d$. 
\[  \makebox {Let  } p_0(x_0) = \sup \{ p: X \in  L^{p}_{loc}  (\RR^d)  \mbox{ in a neighborhood of   $x_0$}  \} . \]   
The function $ p \rightarrow h_{p}(x_0)$ is defined  on  $[1, p_0(x_0) )$ (and possibly also at   $p_0(x_0) $), 
and possesses the following properties: 
\begin{enumerate}
\item It takes values in $\left[-\frac{d}{p_0(x_0)}, \infty \right]$
\item\label{it:dec} It is a decreasing function of $p$. 
% and takes values in $[-d/p, \infty ]$. 
\item The function  $ r \mapsto h_{1/r}(x_0)$ defined on $[0,1]$ is concave. 
\end{enumerate}
Furthermore, Conditions 1 to 3 are optimal, i.e. if  $p_0 \in (1, \infty )$ and $\phi $ is a function defined on $[1, p_0]$ and satisfying
  the above  conditions, 
then there exists $X \in L^{p_0} (\RR^d ) $ such that 
\BE \label{phip} \forall  p \in [1, p_0], \qquad h_{p}(x_0) = \phi (p) . \EE

\ET 

\noindent The proof of the first part of  Theorem~\ref{th:tha} is detailed in \ref{sec:prooftha}. 
The proof of optimality  is detailed in \ref{sec:proofprop2} and relies on the explicit construction of a function $X$ such that the function  $ p \rightarrow h_{p}(x_0)$ can be any function satisfying Conditions 1 to 3 from Theorem \ref{th:tha}. 
This explicit construction also serves as a reference example of singularities where $p$-exponents differ  at $x_0 = 0$ (i.e., the function  $h_{p}(0)$ is not a constant), see Section \ref{sec:skinny} below. \\

\noindent Theorem~\ref{th:tha} calls for the following remarks: 

First, as already stated above, $p$-exponents can by definition take values down to $-d/p$ and hence allow to formalize the notion of {\em negative regularity exponents}.

Second, $p$-exponents for different values of $p$ do in general not  coincide. 
In particular, the H\"older exponent ($p = \infty$) does not in general coincide with  $p$-exponents for  $p < \infty$. 

 Third, as a consequence of the concavity of the function $ r \mapsto h_{1/r}(x_0)$, it follows that, when it takes finite values, $h_{p}(x_0)$ as a function of $p$ is  continuous, except perhaps at endpoints $1$ and $p_0$ (indeed a concave function is continuous except perhaps at the end-points of its domain of definition).
 
Finally, the more difficult problem of understanding which properties are satisfied by the functions of two variables $(x_0, p) \rightarrow h_{p}(x_0)$ remains largely open.  

\subsection{Pedagogical and reference examples}
\label{sec:ex}

Let us now illustrate $p$-exponents on pedagogical examples. 

\subsubsection{Example 1:  Cusp} 

The very reference for local singularity consists of the \emph{cusp} function:
\BE 
\label{cusp} \mbox{If} \quad  \al \notin 2 \NN \qquad
{ \mathcal C}_{\al} (x) = | x-x_0|^\al 
\EE  
It is straightforward to show that the $p$-exponent of ${ \mathcal C}_{\al} $ at $x_0$ does not depend on $p$ and  takes  the constant value $\al$, as is the case with H\"older exponent, i.e., 
$$ h_p(x_0) \equiv \al,  $$
 for all $p \leq p_0$ with 
\begin{equation}
p_0   = 
\begin{cases}
+ \infty & \makebox{ when } \al > 0, \\
  -d/\al  & \makebox{ when } \al < 0. 
\end{cases}
\end{equation}
Examples of cusp singularities with positive and negative regularity are illustrated in Fig.~\ref{fig:cusp}.

\subsubsection{Example 2:  Chirp} 
\label{sec:chirp}

Chirp functions are another classical and commonly studied example of singular behavior. They are defined as
 \BE 
 \label{chirp} 
 \mathcal{C}_{\al,\beta} (x) = | x-x_0|^\alpha \sin \left( \frac{1}{| x-x_0|^\beta}  \right)
 \EE
and serve as the reference for the class of \emph{oscillating functions}.
The application of (\ref{equ:tpe}) immediately yields that 
 $$ h_p(x_0) \equiv \al,  $$
 for all $p\leq p_0$ with
\begin{equation}
p_0  =
\begin{cases}
+ \infty & \makebox{ when } \al \geq 0, \\
-d/\al  & \makebox{ when } \al < 0. 
\end{cases}
\end{equation}
Two chirps with positive and negative exponents are illustrated and analyzed in Fig. \ref{fig:chirp}  (cf. Section~\ref{sec:osc} for further discussions).\\

\subsubsection{Example 3:  Lacunary comb} 
\label{sec:skinny}

For cusp or chirp singularities, $p$-exponents  do not vary with $p$. 
To illustrate the importance of the possibility of making use of a collection of different $p$-exponents, let us now introduce  constructive examples of  univariate functions $F_{ \al, \gamma}: \RR \rightarrow \RR$, referred to as lacunary combs, which highlight a major benefit of using $p$-exponents: 
They enable the refined characterization of certain types of singularities that cannot be revealed when using the H\"older exponent. These functions constitute a key-construction in the proof of Theorem~\ref{th:tha}. 

Let  $\alpha \in \RR  $ and $\gamma  >1$. 
The function $F_{ \al, \gamma}$, defined as:
\BE
\label{equ:Fga}
F_{ \al, \gamma}=
\begin{cases}
2^{-\al j}&\mbox{ for} \; x \in [2^{-j} , 2^{-j} + 2^{-\gamma j}]   \mbox{ for }  \; j\in\NN^+,\\
0 &\mbox{ otherwise},
\end{cases} \EE 
has support on a set $U_\gamma = \bigcup_{j \geq 0}  [2^{-j} , 2^{-j} + 2^{-\gamma j}]$.
Examples of functions $F_{ \al, \gamma} $ with $\alpha>0$ and $\alpha<0$ are plotted and analyzed in Fig. \ref{fig:Fag1}.

The goal is now to derive the function $p \rightarrow h_p (0)$ for $F_{ \al, \gamma}$ at location $x_0 =0$ and to compare it with the H\"older exponent.

When $\alpha \geq 0$, $F_{ \al, \gamma}$ is bounded in  the neighborhood of $x_0 =0$, so that it belongs locally to all $L^p$ spaces with $p\in[1,+\infty]$. 
The function $F_{ \al, \gamma} $ is continuous at 0 where it vanishes (see \ref{sec:proofprop2}).
Moreover, $| F_{ \al, \gamma} (x) |  \leq  | x|^\al$, and this estimate is an equality at points $2^{-j}$. 
Therefore, the H\"older exponent reads $h(0)=\alpha$.

When $\alpha<0$, let $a>0$ and $J$ such that $2^{-J} + 2^{-\gamma J} \leq a \leq 2^{-J+1} ;$ 
then 
\begin{equation}
\label{diverspec1} \frac{1}{a} \int_{-a}^a | F_{ \al, \gamma} (x) |^p dx  = \frac{1}{a}   \sum_{j = J}^\infty 2^{-\gamma j} 2^{-\al p j}  .
\end{equation}
For any $p<-\gamma/\al$, the sum \eqref{diverspec1} converges and therefore $F_{ \al, \gamma}$ belongs locally to $L^p$.
The function $F_{ \al, \gamma} $ is no longer bounded near $x_0 = 0$, and the H\"older exponent at $x_0 = 0$ is thus no longer defined, yet $p$-exponents are well-defined as long as $p< -\gamma/\al$. Note that the minimal regularity condition $f \in L^1_{loc}$ requires that $\al > -\gamma$, which we assume from now on.  

Let us  compute the $p$-exponents of $F_{ \al, \gamma} $ at $x_0 =0$, assuming that $p < -\gamma / \al$ if $\al <0$.
The relation (\ref{diverspec1}), with $ 2^{-J} + 2^{-\gamma J} \leq a \leq 2^{-J+1}$, 
implies that 
\[ \frac{1}{a} \int_{-a}^a | F_{ \al, \gamma} (x) |^p dx  = \frac{1}{a}   \sum_{j = J}^\infty 2^{-\gamma j} 2^{-\al p j}  \sim 2^{-J (\gamma + \al p  -1)} . \] 
This yields
\BE \label{tpx}T^{(p)}_X(a,x_0) \sim a^{ \al + (\gamma -1)/p}. \EE 
If we now take $a$ in the interval $[2^{-j} , 2^{-j} + 2^{-\gamma j}] $, we obtain a quantity which is bounded by the value of $T^{(p)}_X(a,x_0) $ in \eqref{pexpa} at $2^{-j} + 2^{-\gamma j}$, so that 
(\ref{tpx}) still holds.   
We have therefore obtained that, $\forall \al > -\gamma$,  the $p$-exponent of $F_{ \al, \gamma} $ at the origin is   
\BE  \label{pexpfab} \forall p\in [1, p_0) \qquad h_p (0) = \al + \frac{\gamma -1}{p}  \EE  
 with
\begin{equation}
p_0  = 
\begin{cases}
+ \infty & \makebox{ when } \al \geq 0, \\
-\gamma/\al  & \makebox{ when } \al < 0. 
\end{cases}
\end{equation}
The theoretical values and practical estimates of $h_p(0)$ are illustrated in Fig.~\ref{fig:Fag1}.

This example shows that the $p$-exponents of  $ F_{ \al, \gamma} $ differ at $x_0 = 0$ for all values of $p\in [1, p_0 ) $ and thus supply a much richer characterization tool, even in the case when $\al >0$ and the H\"older exponent can be computed.
This will be further discussed in Section \ref{sec:osc}.

Here, we have constructed a function such that  the $p$-exponent is an affine function of $1/p$; the general case (pertaining to the proof to the second part of Theorem \ref{th:tha})  is detailed  in  \ref{sec:proofprop2}.

\section{Wavelet characterization}
\label{sec:wavtool}

\subsection{Discrete wavelet coefficients and $p$-leaders}

It has long been recognized that wavelet coefficients constitute ideal quantities to study the regularity of data (see e.g., \cite{Jaf91, muzyetal94,jmf2,riedi03,Jaffard2004}). 
The characterization of $p$-exponents proposed here relies on the use of the $d-$dimensional discrete wavelet transform (dDWT), whose definition and properties are briefly recalled. \\

\noindent {\bf Mother wavelets.} \quad Let $\{\psi^{(i)}(x)\}_{ i = 1, \cdots 2^d-1} $ denote a family of  oscillating functions with fast  decay and strong joint time-frequency localization properties, referred to as the {\em mother wavelets}.  
Let us further assume that these functions are chosen such that the collection of templates of $\psi^{(i)}$, \BE
\label{wavbasbad}  
2^{-dj/2}\psi^{(i)}(2^{-j}x-k), \hspace{5mm} \mbox{ for} \hspace{3mm} i = 1, \cdots 2^d -1, \hspace{3mm}   j\in \ZZ \hspace{3mm} \mbox{ and} \hspace{3mm}
k\in\ZZ^d 
\EE
dilated to scales $a = 2^{j}$ and translated to space positions $2^{j}k$,
form an orthonormal basis of $L^2 (\RR^d)$ \cite{Mey90I}.
The  $d$-dimensional orthonormal wavelet bases mostly used in practice are obtained by tensor product of univariate orthonormal wavelet basis \cite{cohen}. \\

\noindent {\bf Discrete Wavelet Transform.} \quad  
 The coefficients  of the dDWT of $X$ are defined as
\begin{equation}
\label{eq-wc}
c^{(i)}_{ j,k} =    \displaystyle\int_{\R^d} X(x)  \;  2^{-dj}\psi^{(i)}(2^{-j}x-k) \, dx.
\end{equation}
Note the use of an  $L^1$ normalization for the wavelet coefficients that better fits local regularity analysis and yields the correct self-similarity exponent of the wavelet coefficients for self-similar functions, see \cite{muzyetal94,Jaffard2004,Wendt:2007:E,Abel}. 
For further details on wavelet bases and wavelet transforms, the reader is referred to e.g., \cite{Mallat1998}.\\

\noindent {\bf Uniform regularity and number of vanishing moments.} \quad 
The mother wavelets $\{\psi^{(i)}(x)\}_{ i = 1, \cdots 2^d-1} $ are further required to possess additional regularity and localization properties: 
They are assumed to belong to $C^{ r_\psi } (\RR^d)$ with $r_\psi$ as large as possible.
 When $r_\psi\geq 1$, all their partial derivatives of order at most  $r_\psi$ have fast decay.
Also, the number of vanishing moments $N_\psi $ is defined as a positive integer such that for any polynomial $P$ of degree strictly smaller than $N_\psi $, 
\begin{equation}
\label{eq-vm}
  \int_\R P(x)  \psi^{(i)}(x) dx =  0.
\end{equation}
Both the regularity and the vanishing moment assumptions are required in order to obtain accurate wavelet characterizations  of pointwise regularity:
Let $h_{\max} $ denote the largest smoothness order found in $X$, then a sufficient condition for choosing the mother wavelet reads: 
\begin{equation}
\label{equ:regu}
h_{\max}  < \makebox{ min }(r_\psi, N_\psi). 
\end{equation}
In general, one does not have information concerning a priori regularity of the data, and therefore, one does not know how smooth the analyzing wavelets should be. 
In practice, a rule of thumb consists in using smoother and smoother wavelets, until the outcome no longer depends on the wavelet used, which is interpreted as meaning that sufficient regularity has been reached.  
This can afterwards be confirmed using multifractal analysis tools (see \cite{PART2}), which allow to determine the maximum regularity exponent present in the data.
Further, with orthornormal wavelet bases (such as the so-called "Daubechies"  compactly supported wavelets  \cite{Daubechies}, widely used in applications, and used in this contribution as well as in the companion article \cite{PART2}), one necessarily has $N_\psi  \geq r_\psi$. 
A sufficient (and conservative) condition for accurate wavelet characterizations of pointwise regularity simplifies to $ h_{\max}  < r_\psi $. 

Another (related) consequence of the vanishing moment requirement of major practical importance consists of the fact that wavelets are orthogonal to polynomials: 
The wavelet characterization of either H\"older or $T^p_\al (x_0)$ regularity exponents thus avoids the explicit estimation of the Taylor polynomial. \\

\noindent {\bf\boldmath $p$-leaders.} \quad The multiscale quantities suited to the characterization of $p$-exponents are referred to as  {$p$-leaders} and constructed as the $\ell^p$ norm of a subset of wavelet coefficients $c_\lambda$.
Let $\lambda_{j,k} $, $ k = (k_1, \ldots, k_d)$ and $j \in \ZZ$, denote a dyadic cube 
\[   \lambda_{j,k} \; : = \left[ 2^j k_1,  2^j (k_1+1)  \right) \times \cdots \times  \left[ 2^j k_d,  2^j (k_d+1) \right). \]
Let $C\lambda$ denote  the cube homothetical to the dyadic cube $\la$, with same center,  but $C$ times wider.
Accordingly, $c_{\lambda} \equiv c_{ j,k}$, and $\psi^{i}_{\lambda} \equiv \psi^{(i)}_{ j,k}$.\\

\begin{figure}[tb]
\centering
\setlength{\tabcolsep}{0pt}
\includegraphics[width=0.45\linewidth]{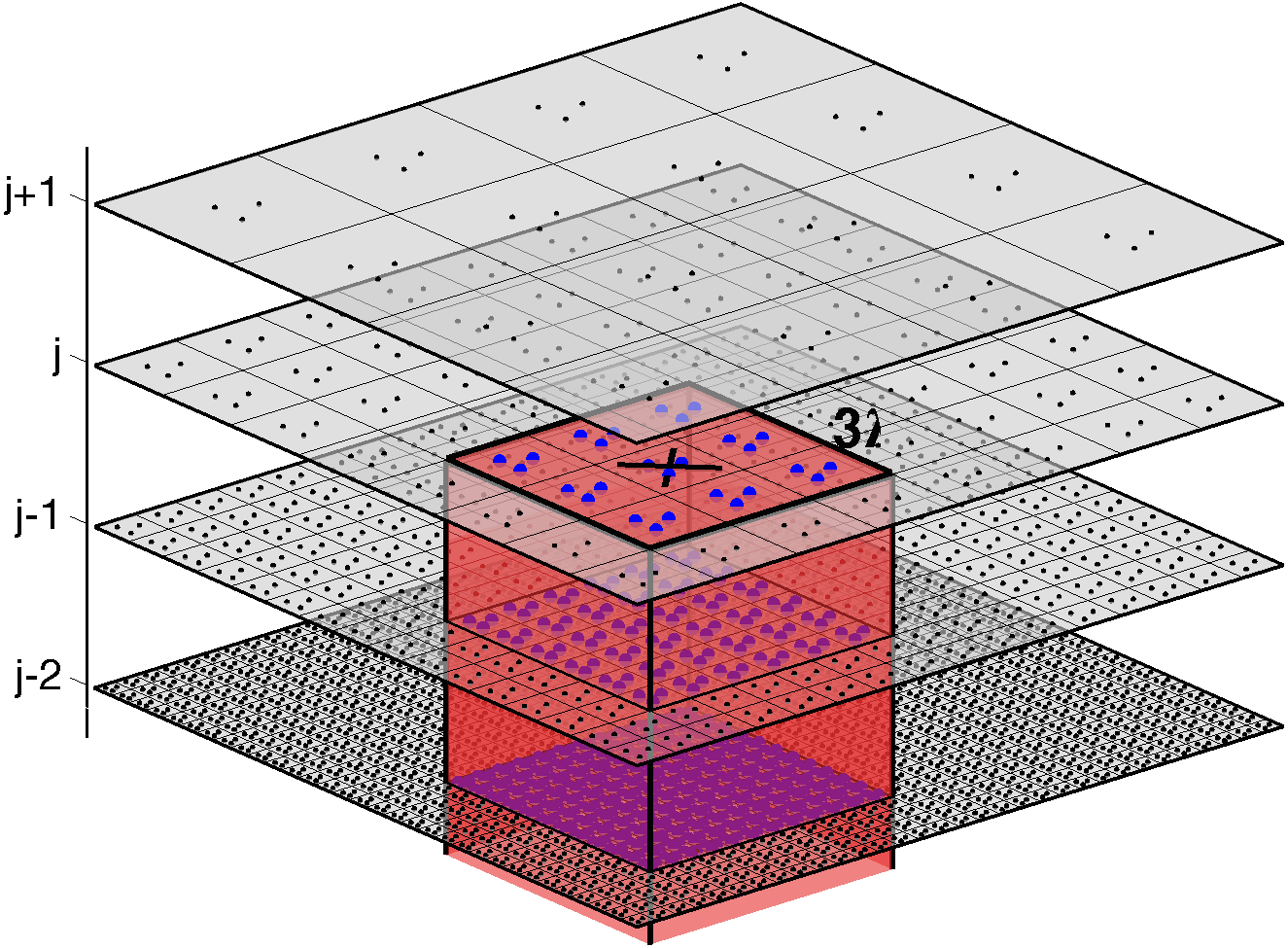}\vspace{-3mm}
\caption{\label{fig:leaders} {\bf\boldmath Definition of $p$-leaders.}}
\end{figure}

 If $X \in L^p_{loc} (\RR^d)$,  the $p$-leaders are defined as: 
\begin{equation} 
\label{pleaders} 
 \pL(j,k)  \equiv \pL_\la= \left(   \sum_{j' \leq j, \;  \lambda' \subset 3 \lambda} \sum_{i=1}^{2^d-1} \big| c^{(i)}_{ j,k} \big|^p \, \, 2^{-d(j-j')} \right)^{1/p}
\end{equation}
where $j'\leq j$ is the scale associated with the sub-cube $\la'$ included in  $3 \lambda$  (i.e. $\la'$ has width $2^{j'}$). 
The $\ell^p$ norm is thus taken on all cubes $\la'$  of scale at most $2^{j}$, which are included either in $\la$ or in its $3^d -1$ closest neighbors, as illustrated in Fig. \ref{fig:leaders}.

When $p = +\infty $ (and $X \in  L^{\infty}_{loc} (\RR^d)$), {$p$-leaders} boil down to the classical wavelet leaders
$$\ell_\lambda = \sup_{i \in [ 1,\ldots, 2^d-1], \;   j' \leq j , \; \lambda ' \subset 3 \lambda}  |c_{\lambda'}^{(i)}|$$
used to characterize H\"older exponents \cite{Jaffard2004,Wendt:2007:E}.

\noindent{\bf Remark.\quad} Note that an equivalent definition of $p$-leaders is given by
\begin{equation} 
\label{pleaderssup} 
 \pL(j,k) = \left(   \sum_{j' \leq j, \;  \lambda' \subset 3 \lambda} d_{j',k'}^p \, \, 2^{-d(j-j')} \right)^{1/p}
\end{equation}
where $d_{ j,k} = \sup_{i = 1, \cdots 2^d -1} \left| c^{(i)}_{ j,k} \right|$.

\subsection{Global regularity and wavelet coefficients} 

\noindent {\bf Local boundedness.} \quad As already discussed, the H\"older exponent suffers from a well-known and documented limitation: 
It can be used as a measurement of pointwise regularity for locally bounded functions only, see \cite{JAFFARD:2010:A}.  
Whether empirical data can be well-modeled by locally bounded functions, or not, can be practically determined through the computation of the {\em uniform H\"older exponent}  $\Hmin$, using the following wavelet characterization:
 \BE 
\label{caracbeswav3hol}  
\Hmin =\liminf_{j \rightarrow - \infty} \;\;  
 \;\;  \frac{ \log 
\left(  \displaystyle  \sup_{ i,k}  \big| c^{(i)}_{ j,k} \big|     \right) }{\log (2^{j})}. 
\EE
Indeed, if {  $\Hmin > 0$}, then $X$ is  a continuous  function,  whereas, if  {  $\Hmin < 0$}, then $X\notin L^\infty_{loc}$,  see \cite{MandMemor,Bergou}. 
For numerous real world applications, the restriction $\Hmin > 0$  constitutes a severe limitation, cf. \cite{MandMemor}.
From a practical point of view, mother wavelets satisfying $r_\psi>\Hmin$ are required for an accurate  estimation of $\Hmin$.\\

\noindent{\boldmath \bf Local $L^p$ regularity.} \quad Wavelet coefficients can furthermore be used to assess whether $X$ locally belongs to $L^p$ or not, i.e., to verify the condition that is a priori required for using the corresponding $p$-exponent \cite{Bergou,Abel,MandMemor}.
Let $S_c(j,p)$ denote the wavelet structure function, defined as the space/time averages of the magnitude of wavelet coefficients raised to a positive power $p  >0 $:
\begin{equation}
\label{equ-WSF}
S_c(j,p) = 2^{dj} \displaystyle\sum_{  k }\sum_{i=1}^{2^d-1}  \big| c^{(i)}_{ j,k} \big|^p. 
\end{equation}
Let  $ \eta_X (p) $ denote the {\em wavelet scaling function}, defined as
\BE \label{defscalond}  \forall p >0, \hspace{6mm} 
 \eta_X (p) =   \displaystyle \liminf_{j \rightarrow - \infty} \;\; \frac{\log \left( S_c(j,p)  \right) }{\log (2^{j})} \EE 
Because of its function space interpretation in terms of Besov spaces (see Section \ref{sec:ppluspetun} below), $ \eta_X (p)$ does not depend on the wavelet basis as long as $r_\psi>  | \Hmin | $ and $p\geq 1$. 
Then when $p\geq 1$: 
\begin{equation}
\label{equ:etap}
\left. \begin{array}{rl}  \mbox{ If} \; \eta_X(p) >0 & \mbox{then} \; X \in L^p_{loc}, \\ 
\mbox{ if} \; \eta_X (p) < 0 & \mbox{then} \; X \notin L^p_{loc} . \\ 
\end{array}\right\}
\end{equation}

The wavelet scaling function $\eta_X (p)$ must therefore be computed prior to applying $p$-exponent characterization to determine the range of values of $p$ suitable for analysis.

\subsection{$p$-exponent regularity characterization with $p$-leaders}\label{subsec:pcharac}

Let $p > 1$. 
When $\eta_X (p) >0$ and $r_\psi > h_p (x_0)$, the $p$-exponent $h_p(x_0)$ defined in \eqref{equ-pexp} can be recovered from $p$-leaders \cite{JaffToul, Bergou,JaffMel}: 
\begin{equation} 
\label{carachqf}  
h_p(x_0)  =  \liminf_{j \rightarrow - \infty}  \frac{ \log  \left(  \pL_{\lambda_{j,k} (x_0)} \right)  }{\log (2^{j})}. 
\end{equation}
Note that  the characterization provided by (\ref{carachqf}) does not require the computation of the Taylor polynomial, and that it
extends to any $p$ the previously obtained wavelet leader characterization of H\"older exponents (valid when $\Hmin >0 $) \cite{Jaffard2004,Wendt:2007:E}:
\[ h(x_0) =\liminf_{j \rightarrow - \infty} \frac{ \log \left(    \ell_{\lambda_{j,k} (x_0)}  \right) }{\log (2^{j})}.\]  

\subsection{ The case $0< p <1$}  
\label{sec:ppluspetun}

\noindent {\bf\boldmath Limitation $p \geq1$.} \quad The restriction $p \geq 1$ may constitute a severe drawback in applications as it implies that pointwise singularities with regularity smaller than $-d$ cannot be dealt with. 
As already pointed out, the restriction $p \geq 1$ stems from $L^p$ spaces being ill-defined when $0<p<1$.
In particular, even when the wavelet basis belongs to the Schwartz class and $X$ satisfies $\int |X(x) |^p \; dx < \infty$ for $p<1$, its wavelet coefficients might not be well-defined: 
Indeed for example, for the function $X(x) = |x|^{-2}$, which belongs e.g. to $L^{1/3}$, the scalar product $\int X\varphi$ is in general undefined even if $\varphi$ is $C^\infty$, unless $\varphi(0)=\varphi'(0)=0$.  \\

\noindent {\bf\boldmath Replacing $L^p$ spaces with Besov spaces.} \quad However, this issue has been dealt with mathematically by the substitution of $L^p$ spaces by Besov spaces, see \cite{cohen,Mey90I}. 
We sketch below only the basic notions together with practical implications that are relevant to pointwise regularity. 

Let $p >0$ and $s \in \RR$. 
A function $X$ belongs to the  Besov space $B^{s,\infty}_p$ if its wavelet coefficients satisfy
\[ \exists C\;  \forall j \qquad 2^{dj} \sum_{ k} ( d_\la )^p \leq C 2^{spj}  \]
(where the sum bears on all dyadic cubes of width $2^j$). 
For applications, it is important to note that the original requirement that the mother wavelet belongs to the Schwartz class can be relaxed to wavelets with sufficient regularity and vanishing moments: 
$r_\psi>|s|\mbox{ if }p\geq 1$ and $r_\psi >  s  >d(2/p-1)-r_\psi\mbox{ otherwise}$ \cite{Bourdaud1995}. 

When $p \geq 1$, Besov spaces are closely related with $L^p$ spaces via  the following embeddings:
\[ \forall \ep >0, \qquad B^{\ep ,\infty}_p \hookrightarrow L^p \hookrightarrow B^{-\ep ,\infty}_p. \] 
In particular, the condition $\eta_X (p) >0$, used as a practical sharp condition implying that data locally belong to $L^p$, is equivalent to 
requiring that  $X $ locally belongs to $ B^{\ep ,\infty}_p$ for an $\ep >0$. 
The requirement $X \in L^p$ can thus  systematically be replaced by $X \in B^{0,\infty}_p$ in all results reported so far, with no change in any of them, yet with the
 benefit that it readily extends down to $ p>0$.  
The definition of $p$-leaders given by (\ref{pleaders}) remains unchanged, and $p$-exponents for $0 < p <1$ can now directly be defined through (\ref{carachqf}).  \\

\noindent {\bf Admissible distribution.} \quad From now on, we therefore assume that there exists $p_0 >0$ such that $X$ locally belongs to  $ B^{\ep , \infty}_{p_0}$ for an $\ep >0$ (i.e., $\eta_X(p_0)>0$), so that $p$-exponents can be defined, at least,  for $0< p \leq p_0$.
When this condition is fulfilled, $X$ is said to be an \emph{admissible tempered distribution} (see Section \ref{sec:wgn} for an example of a random distribution which is not admissible). The extension of 
Theorem \ref{th:tha}, to the full range  $p\in (0, + \infty]$, will be the subject of a forthcoming paper. \\

\noindent{\boldmath \bf Chirps with negative regularity.\quad} To illustrate the relevance of using $p$-exponents with $0<p<1$, let us consider the example of chirp functions (cf. Section~\ref{sec:chirp}).
When $\al \leq -1$,  $\mathcal{C}_{\al,\beta} $ does not belong to $L^1_{loc}$.
However, ${ \mathcal C}_{ \al,\beta}$ locally belongs to $B^{0,\infty}_p$ when $p \leq -1/\al$ and $h_p(x_0)=\al$. 

\subsection{Illustrations, examples and counter-examples}

\subsubsection{Illustrations and examples} 

\begin{figure}[htp]
\centering
\includegraphics[width=\linewidth]{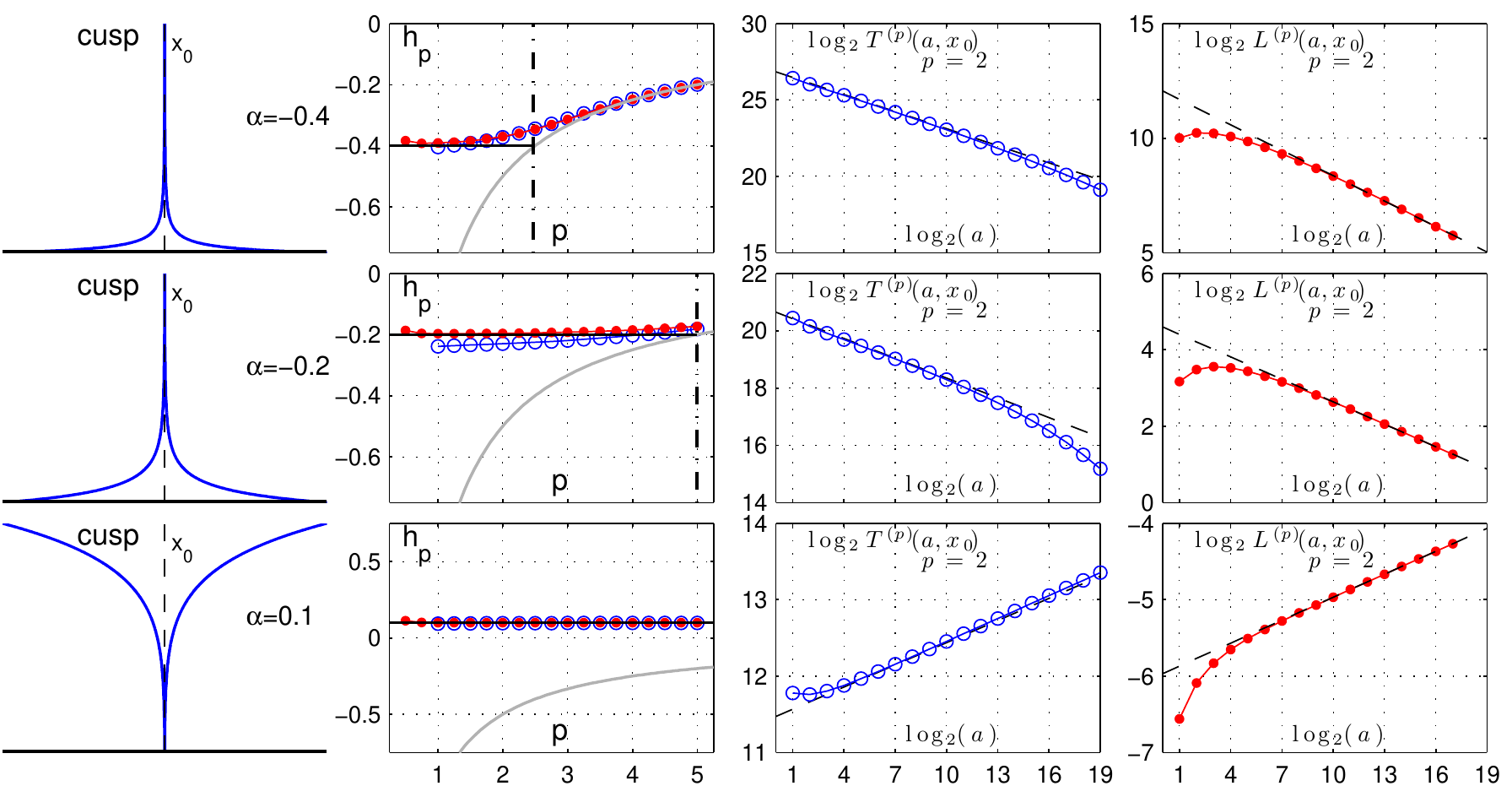}
\caption{\label{fig:cusp} {\bf\boldmath $p$-exponents of cusps.} Left column: Cusps with positive and negative $p$-exponent. Second column: Estimates of $h_p(x_0)$ using $T_\alpha^{(p)}$ definition (\ref{pexpa}-\ref{equ-pexp}) (blue circles) and $p$-leaders (red discs), theoretical exponent (black solid) and $L^p$ limit (grey). Third and fourth columns: $\log_2T_\alpha^{(p)}(a,x_0)$ and $\log_2\pL(a,x_0)$ as a function of $\log_2(a)$.}
\end{figure}

\begin{figure}[htp]
\centering
\includegraphics[width=\linewidth]{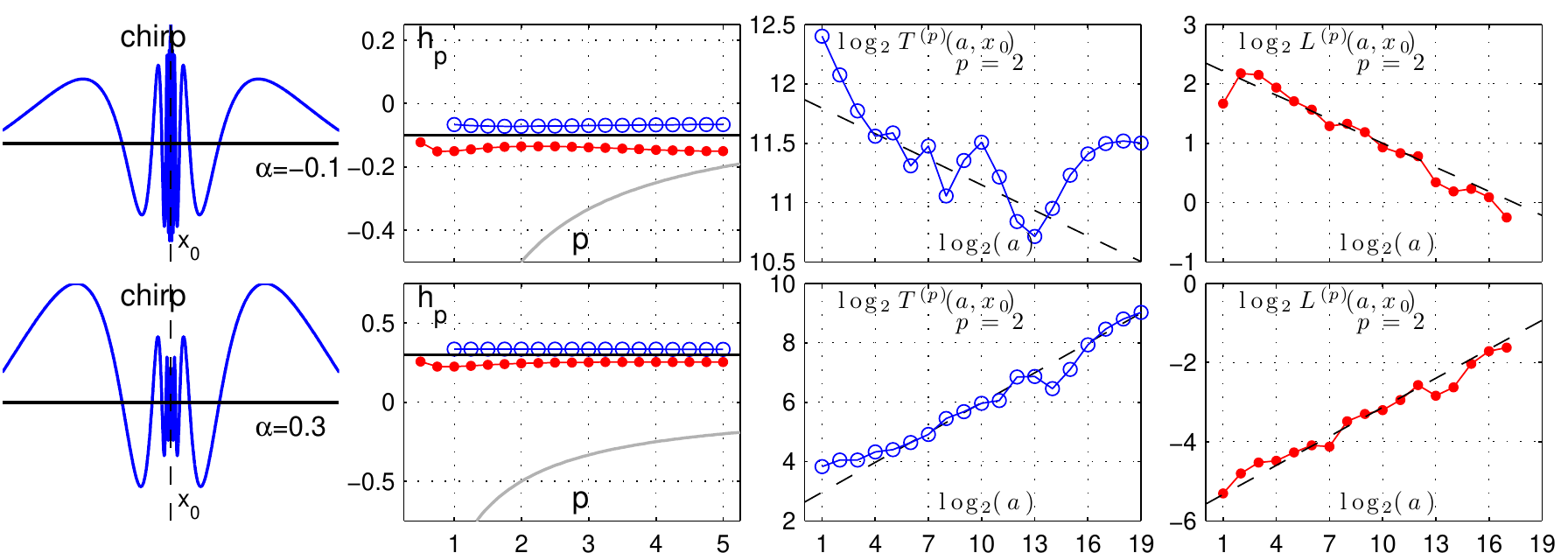}
\caption{\label{fig:chirp} {\bf\boldmath $p$-exponents of chirps.} Left column: Chirps with positive and negative $p$-exponent. Second column: Estimates of $h_p(x_0)$ using $T_\alpha^{(p)}$ definition (\ref{pexpa}-\ref{equ-pexp}) (blue circles) and $p$-leaders (red discs), theoretical exponent (black solid) and $L^p$ limit (grey). Third and fourth columns: $\log_2T_\alpha^{(p)}(a,x_0)$ and $\log_2\pL(a,x_0)$ as a function of $\log_2 a$.}
\end{figure}

\begin{figure}[htp]
\includegraphics[width=\linewidth]{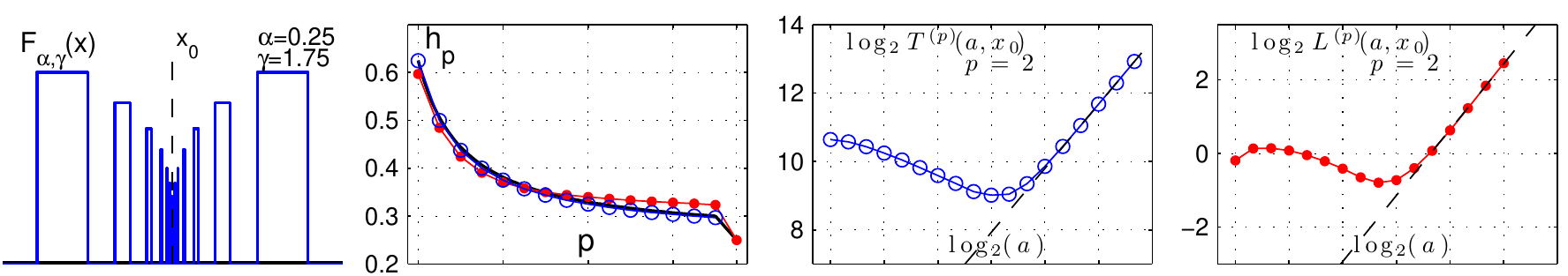}\\
\includegraphics[width=\linewidth]{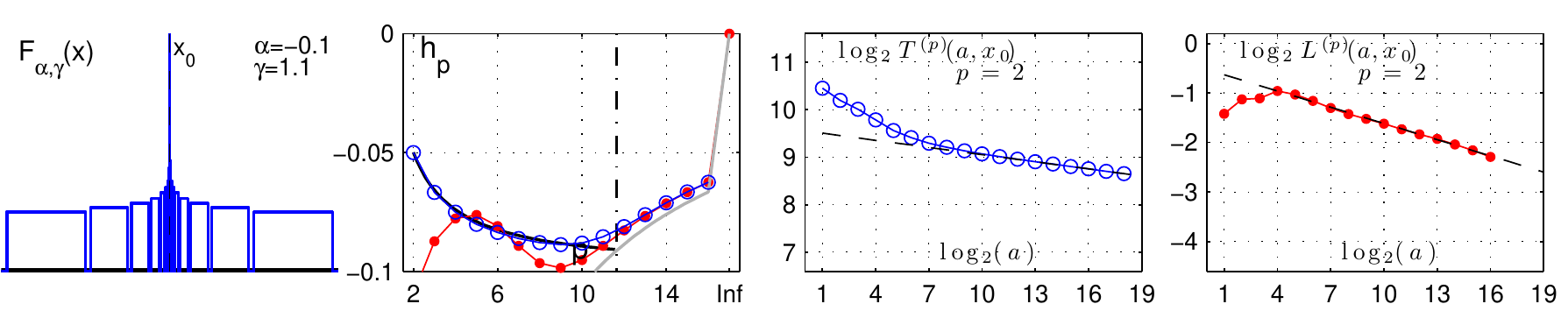}
\caption{\label{fig:Fag1} {\bf\boldmath Lacunary comb $F_{\alpha,\gamma}(x)$.} 
Two functions $F_{\alpha,\gamma}(x)$ with positive (top row) and negative (bottom row) $p$-exponents (left column).
Second column: Estimates of $h_p(x_0)$ using $T_\alpha^{(p)}$ definition (\ref{pexpa}-\ref{equ-pexp}) (blue circles) and $p$-leaders (red discs). Third and fourth columns: $\log_2T_\alpha^{(p)}(a,x_0)$ and $\log_2\pL(a,x_0)$ as a function of $\log_2 a$.
The black solid line indicates the theoretical $p$-exponent, the grey solid line the limit for the singularity to be in $L^p_{loc}$.
}
\end{figure}

The power-law behaviors underlying the wavelet characterization of $p$-exponents (cf. \eqref{carachqf}), and thus their $p$-leader based estimation, are illustrated on the reference examples defined in Section~\ref{sec:ex}, and compared with an estimation relying on a direct use of the definition of the $T^p_\alpha $ regularity (cf. Section~\ref{sec:pexpodef} and \eqref{pexpa} and \eqref{equ-pexp}). 
Fig. \ref{fig:cusp} shows cusp singularities (cf. \eqref{cusp}) with positive and negative exponents.
Fig. \ref{fig:chirp}  illustrates the analysis of chirps (cf. \eqref{chirp}), with positive and negative exponents.
Fig. \ref{fig:Fag1} addresses the analysis of lacunary combs  $F_{\al,\gamma}$ (cf. \eqref{equ:Fga}), with positive and negative exponents.

\subsubsection{A counter example: white Gaussian noise} 

\label{sec:wgn}

Despite its being very broad,  the class of admissible distributions does not contain all tempered distributions. 
Notably, the practically very natural example of  univariate independent identically distributed (IID) Gaussian random variables, also referred to as white Gaussian noise (wGn), satisfies $\Hmin =  -1/2$.
However, its wavelet coefficients read 
\[  c_{j,k} = 2^{-j/2} \chi_{j,k} , \]
where the $\chi_{j,k}$ are IID zero-mean Gaussian random variables.
This implies that the wavelet scaling function reads $\eta_X (p) = -p/2 < 0$, for all $p >0$, and therefore that $p$-exponents are not fitted to analyze WGN.  
The same holds for fractional Gaussian noise (fGn), the increment process of fractional Brownian motion, the only Gaussian self-similar process. 
Whatever the value of the self-similarity parameter $0<H<1$, $\eta_X (p) = -pH < 0$ for all $p >0$.

\section{Pointwise singularity classification}
\label{sec:pExp}  

\subsection{Motivation}

Traditionally, singular behaviors are categorized using two classes, non oscillating versus oscillating singularities, cusps and chirps (as defined in Section~\ref{sec:ex}) constituting the reference models for each class respectively. 
Deciding whether real-world data contain oscillating singularities or not is of practical importance, as it may change the analysis of the underlying physical and biological mechanisms at work behind data \cite{GRETSI15-pL,Bandt2015}. 
This is notably the case in hydrodynamic turbulence modeling, where the presence/absence of oscillating singularities may permit to validate/falsify various vortex stretching mechanisms \cite{Frisch1995}.
Such a classification currently relies on the use of two exponents, the H\"older and the oscillation exponents, and is best illustrated using the concept of integration (cf. e.g., \cite{Bergou}). 

Let us compare the 1D chirp function ${ \mathcal C}_{\alpha, \beta} $ (cf. \eqref{chirp}) to its primitive ${ \mathcal D}_{\alpha, \beta} $. 
For the former, the $p$-exponents at $ x_0$ read $h_p(x_0)=\alpha$ (for $p\in (0, +\infty ] $ if $\al \geq 0$ and  $p \in (0, -1/\al )$ if $\al <0$) and are hence independent of $\beta$ and $p$.
For the latter, a simple integration by parts yields 
 \[ { \mathcal D}_{\alpha, \beta}  (x) = \frac{| x-x_0|^{\alpha + \beta +1} }{\be} \cos \left( \frac{1}{| x-x_0|^\beta}  \right) + C+ \mathcal{O} \left( | x-x_0|^{\al + 2 \beta+2}  \right),    \]
 thus showing that the increase of any $p$-exponent after integration reads $\be +1$. 
This is in contradistinction with the cusp case (\ref{cusp}), where the increase after integration is exactly $1$.  
This simple computation indicates that cusps and chirps can be discriminated by integration and paves the way towards the definition of an {oscillation exponent} $\be$ (cf. e.g., \cite{Bergou} for the case of H\"older exponents). 

However, both cusps and chirps are characterized by $p$-exponents that do not depend on $p$ and can hence not be considered different with respect to $p$-exponents, but instead depart from the lacunary comb example (cf. Section~\ref{sec:skinny}).
A general interpretation of these two different behaviors is provided in \cite{Bandt2015}.

This naturally leads to define the new class of {\em $p$-invariant} singularities. 
\BD 
\label{def:pi}
\label{th:propo4bis} 
 Let $X: \RR^d \rightarrow \RR$ be an admissible  distribution. 
The singularity at $x_0$ is said to be $p$-invariant if and only if
\BE \label{}  
\forall p \in (0,  p_0(x_0)), \quad h_{p} (x_0)= h. 
\EE
In other words, the function $p \rightarrow h_p (x_0) $ is constant in the interval $(0, p_0(x_0) )$.  
\ED

The remainder of the section aims to consider the extent to which the flexible setting of $p$-exponents permits to enrich the classification of singular behaviors, as illustrated in Fig.~\ref{fig:singularitymap}.

\subsection{Canonical vs. oscillating singularities} 
\label{seccanon}
\subsubsection{Fractional integration} 

While examples of the previous section were 1D functions, the multidimensional setting, $x \in \RR^d$ with $d \geq 2$, is considered, thus requiring to replace the concept of \emph{primitive function} with that of {\em fractional integral}:  
\BD 
\label{def:fi}
Let  $X$ be a tempered distribution defined on $\RR^d$.
The fractional integral of order $s$ of $X$, denoted by  $X^{(-s)}$, is defined as the convolution operator $(Id-\Delta )^{-s/2}$.
Equivalently, in the Fourier domain, it corresponds to the multiplication by the function  $(1+|f|^2)^{-s/2}$.
The $p$-exponent of $X^{(-s)}$ at $x_0$ is well defined on condition that $\eta_X (p)> -s p$ and is denoted by $h_{p,s} (x_0)$.
\ED
The condition $\eta_X(p)> -s p$, sufficient to insure that  $X^{(-s)}$ locally belongs to $L^p$ when $p \geq 1$ (or to $B^{ \ep, \infty}_p$ for an $\ep >0$ when $p < 1$),
 follows from the Besov space interpretation of the wavelet scaling function, see \cite{Mey90I}. 
 
Generalizing \eqref{carachqf}, the exponent $h_{p,s} (x_0)$ can be characterized by 
 \begin{equation} 
\label{carachqft}  
h_{p,s}(x_0)  =  \liminf_{j \rightarrow - \infty}  \frac{ \log  \left(  \pLt_{\lambda_{j,k} (x_0)} \right)  }{\log (2^{j})}, 
\end{equation}
with $(p,s)$-leaders defined as:  
\begin{equation} 
\label{ptleaders} 
\pLt (j,k)  \equiv \pLt_\la= \left(   \sum_{j' \leq j, \;  \lambda' \subset 3 \lambda} \sum_{i=1}^{2^d-1} \big( 2^{sj'} \; \big| c^{(i)}_{\la'} \big| \big)^p \, \, 2^{-d(j-j')} \right)^{1/p}.
\end{equation}
Essentially, this means that, as regards pointwise regularity, taking a fractional integral of order $s$ is equivalent to multiplying the wavelet coefficients by $2^{js}$, see \cite{JaffToul,JAFFARD:2010:A} for a mathematical justification of this heuristic. 

\subsubsection{Canonical singularity} 

The notion of \emph{canonical singularity} is defined as follows: 
\BD  
\label{def:defcansing} 
Let $p$ and $h_p$ be such that $p < p_0(x_0)$ and $h_p > -d / p$. 
An admissible  distribution  $X: \RR^d \rightarrow \RR$ has a canonical singularity of exponent $h_p$ at $x_0$  if  $\eta_X (p) >0$ and if there exists  $s >0$, such that 
\BE
\label{equ-hpg}
h_{p,s} (x_0)= h_p  + s. 
\EE
\ED
A similar notion had previously been considered only in the H\"older case, i.e., when $p= +\infty$, see \cite{me98}. 
The following result shows that canonical singularities are $p$-invariant singularities (thus motivating the choice of the terminology).

\BT 
\label{th:propo4bis} 
Let $X: \RR^d \rightarrow \RR$ be an admissible distribution with  a canonical singularity at $x_0$. 
Then
\BE \label{equ-hpgbis}  \forall p \in (0,  p_0], \quad \forall s \geq 0, \qquad h_{p,s} (x_0)= h_p  + s. 
\EE
In particular, $p \rightarrow h_p (x_0) $ is constant in the interval $(0, p_0(x_0) )$.  
\ET  
The proof of Theorem \ref{th:propo4bis} is detailed in \ref{sec:proofprop4}.  
Note that it does not rely on the assumption that $p_0(x_0) = + \infty$, and its conclusion also applies for $0< p < 1$. \\

\noindent{\bf Example 4: Self-similar  distributions.\quad}
Besides cusps, further examples of canonical singularities are supplied by {\em deterministic self-similar functions or distributions}, which are defined as follows. We only consider the one-variable case here.

\BD Let  $a>1$ and $\al \in \RR$. A nonvanishing tempered distribution $X$ defined on $\RR$ is self-similar of  scaling ratio  $a$ and exponent $\al$ at  $x_0 $ if 
the two distributions $a^\al X( x-x_0 ) $ and $X(a( x-x_0))$ coincide.
\ED 
The following result  shows that self-similar functions  supply simple examples of canonical singularities. 

\BP
\label{prop:ss}
 Let $p \geq 1$ and $X \in L^p_{loc} (\RR )$ be a self-similar function  at $x_0$ of  scaling ratio $a>1$ and exponent $\al$. If $X$ differs from a polynomial in the neighbourhood of $x_0$,  then   $p > -1/\al $,  $X$ has a canonical singularity at $x_0 $ of exponent $\al$, and it can be written under the form 
 \BE \label{eq:prop3}  
  X(x) = 
\left\{ 
\begin{array}{lll}
 |x-x_0|^\al \omega_+ (\log (x -x_0)) &  \mbox{ if} & x >0   \\
  |x-x_0|^\al \omega_- (\log (-(x -x_0)) ) &  \mbox{ if} &  x <0  , \\
\end{array}
\right.
 \EE
where  $ \omega_+   $ and  $ \omega_-  $  are  $\log(a)$ periodic functions in $L^p$. 
\EP 
The proof of Proposition \ref{prop:ss} is given in \ref{proof:selfsimilar}.

Simple examples of self-similar functions or distributions are provided by the cusps ${ \mathcal C}_{\alpha}$ when $\al \notin \ZZ$ (defined by (\ref{cusp})  for $\al > -1$ and in \cite[Section 3.4]{Bandt2015}  for $ \al \leq -1$), in which case  any $a>1$ will fit. 
Note that, if $\al > -1$,   these fall in the case covered by   Proposition~\ref{prop:ss} with  $\omega_+$ and  $\omega_-$ set to  constants.
Several examples of cusps are plotted and analyzed in Fig.  \ref{fig:cusp}. 
Another important example is provided by  singularities with complex exponent   \cite{johanson2001} 
$$ |x-x_0|^{\al+i \beta} = |x-x_0|^{\al} \exp(-i \beta \log |x - x_0|) . $$

Other examples of self-similar distributions which are not functions are supplied by the measures 
\[ \sum_{n \in \ZZ}  a^{n\al} \delta_{a^n}, \qquad \mbox{ for} \quad  a >1 \quad  \mbox{ and } \quad \al > 0. \] 
Note that,  if  $X$ is a self-similar distribution of scaling ratio  $a$ and exponent $\al$, then its derivative $X'$  (in the sense of distributions) is self-similar of scaling ratio  $a$ and exponent $\al -1$; thererefore,  by derivation, we can deduce from the previous examples   self-similar distributions of arbitrary negative order.  

\subsection{Oscillating singularities: Lacunary versus balanced} 
\label{sec:osc}

\subsubsection{Oscillating singularities}

\BD  
\label{def:defoscsing} 
Singularities that are not canonical are called oscillating singularities. 
\ED

The classical examples of oscillating singularities  are provided by the chirps ${ \mathcal C}_{\alpha,\beta}$ in (\ref{chirp}).
Their $p$-exponents are given by $h_p(x_0)=\alpha$ (for any $p\in (0, +\infty ] $ if $\al \geq 0$, and for $p \in (0, -1/\al )$ if $\al <0$). 
The functions ${ \mathcal C}_{\alpha, \beta}$ therefore supply examples of singularities which are not canonical (they are oscillating) yet are $p$-invariant. 

\subsubsection{Lacunary versus balanced singularities}

However, in general, oscillating singularities need not be $p$-invariant and this large class of singularities  can be further refined by the use of $p$-exponents.
Let us consider again the function $F_{ \al, \gamma}$ defined in Section \ref{sec:pexpodef} and  \eqref{equ:Fga}. 
It has an oscillating singularity at $x_0$: 
Indeed, the primitive of $F_{ \al, \gamma}$ has exponent $\al + \gamma$ at $x_0=0$. 
This directly follows from the fact that its integral on the interval $[2^{-l} , 2^{-l} + 2^{-\gamma l}] $ is  $2^{-(\al + \gamma )l} $ and  that $F_{ \al, \gamma}$  vanishes elsewhere.  
The major virtues of the examples  $F_{ \al, \gamma} $, as well as of the more general case developed in \ref{sec:proofprop2}, consists in showing that observing $p$-exponents that differ is deeply tied to the notion of \emph{spatial/temporal lacunarity}: Indeed, the support of the function $F_{ \al, \gamma} $ is a set which is ``scarce'' at $x_0$ (cf. \cite{Bandt2015} for the technical formulation of scarcity). 

This leads to refine the classification of oscillating singularities.
\BD
\label{def:LacOsc} Let $X$ be an admissible distribution, with an oscillating singularity at $x_0$. 
This oscillating singularity is balanced if the function $p \mapsto h_p(x_0)$  is constant in its interval of definition, and is lacunary otherwise.
\ED

\vskip3mm

\begin{figure}[t]
\centering
\setlength{\tabcolsep}{0pt}
\includegraphics[width=0.4\linewidth]{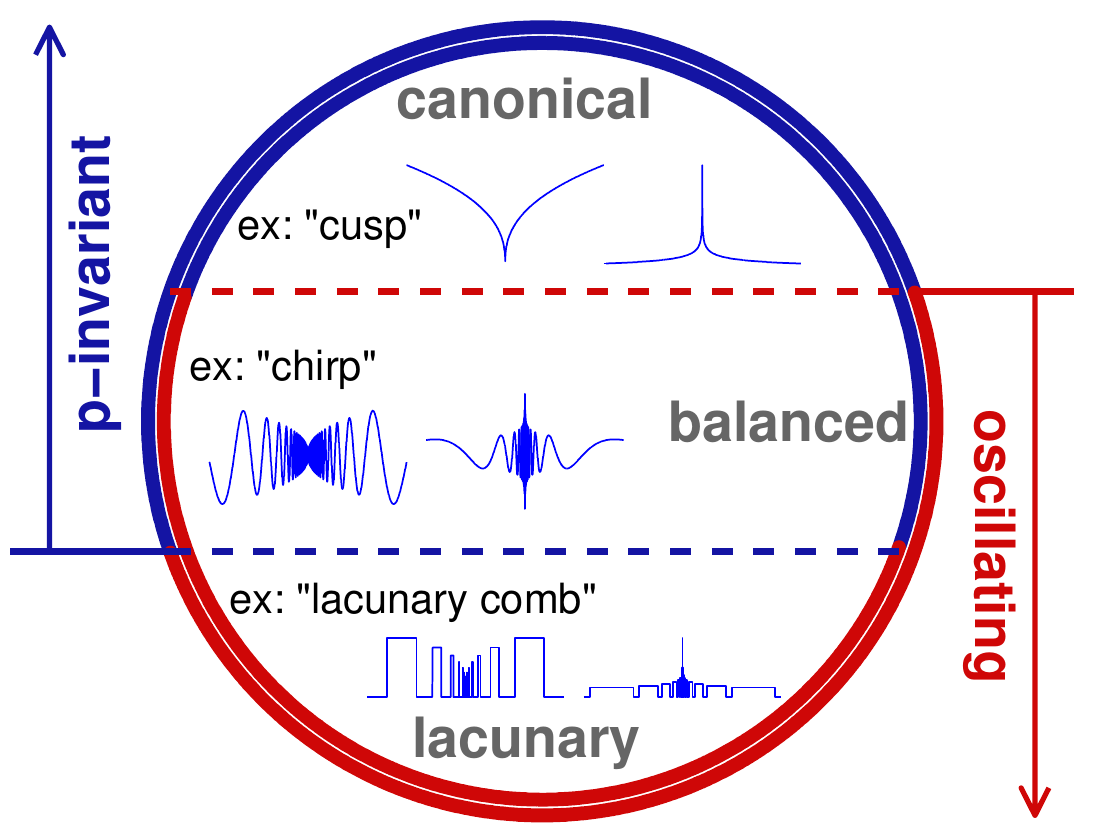}
\caption{\label{fig:singularitymap} {\bf Classification of singularities.}}
\end{figure}

The $p$-exponent based classification of singularities provided by Definitions \ref{th:propo4bis}, \ref{def:defcansing} \ref{def:defoscsing} and \ref{def:LacOsc} is schematically illustrated in Fig. \ref{fig:singularitymap}.

\subsection{Negative regularity}
\label{sec:negative}

\subsubsection{Motivation}

Let us now return to the original motivation of analyzing negative regularity singularities. 
As mentioned in the introduction, the classical strategy to handle negative singularities is to compute $\Hmin $ (using \eqref{caracbeswav3hol}) and, whenever $ \Hmin < 0 $ to perform  a fractional integration of order $s > - \Hmin $. 
Indeed,   the function $X^{(-s)}$ satisfies $ h^{min}_{X^{(-t)}} = h^{min}_{X} + s $, which is positive,  so that the H\"older exponent  can be used  in order to measure its pointwise regularity. 
A major advantage of this approach stems from its not relying from any a priori assumption. 
A significant limitation stems from the fact that (fractional) integration may mask some of the regularities existing in original data. 
Let us examine how the alternative $p$-exponent approach may overcome such limitations.

\subsubsection{Example 5: Cusp plus Chirp}

We now construct a simple toy-example that will show why pointwise regularity information can be lost in (fractional) integration, when data contain oscillating singularities and negative regularity, cf. \cite{ABRY:2011:A}. 
Let 
\[ X_1(x) = |x-x_0|^{\gamma}, \; X_2(x) =  |x-x_0|^{\al } \sin\left( \frac{1}{|x-x_0|^\be} \right), \; \mbox{and} \; X(x) = X_1(x) + X_2(x)  , \] 
with $ \al < \gamma < \frac{\al}{1 + \be } <0$. 
The condition $\al < \gamma $  implies that the dominating singularity at $x_0$  is supplied by $X_2$.
The exponents $\Hmin$ of $X_1$ and $X_2$ are respectively $\gamma$ and $\frac{\al}{1 + \be } $, thus yielding $\Hmin = \gamma $ for $ X $.
Interestingly, $\Hmin$ is thus larger than $\alpha$ and hence no longer yields a lower bound of possible pointwise exponents, as is the case when only singularities with positive regularity  are present. 
A H\"older regularity based analysis of $X$ requires performing a fractional integration of order $ s > -\gamma$. Then, the H\"older exponent of $X_1$ at $x_0$ is shifted to  $ \gamma +s $, and to $\al + s(1+\beta)$  for $X_2$, which is strictly larger than $ \gamma +s $. 
Therefore, the oscillating part $X_2^{(-s)}$  is  dominated by the cusp part $X_1^{(-s)}$, and the information associated to the initially dominant oscillating singularity cannot be recovered. 
Note that in this example, $p$-exponents are equally affected by this limitation caused by the fractional integration.

\begin{figure}[t]
\centering
\begin{tabular}{rr}
\footnotesize sum of cusp and chirp singularities & \footnotesize fractional integral of order $s=0.3$\\
\includegraphics[height=28mm]{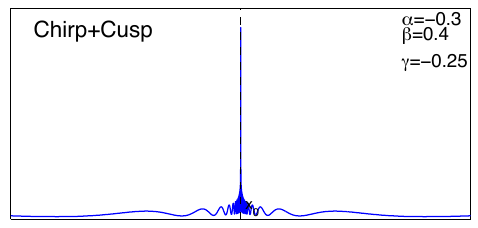}&
\includegraphics[height=28mm]{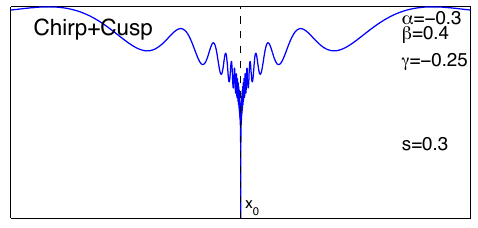}\vspace{-1mm}\\
\includegraphics[height=30mm]{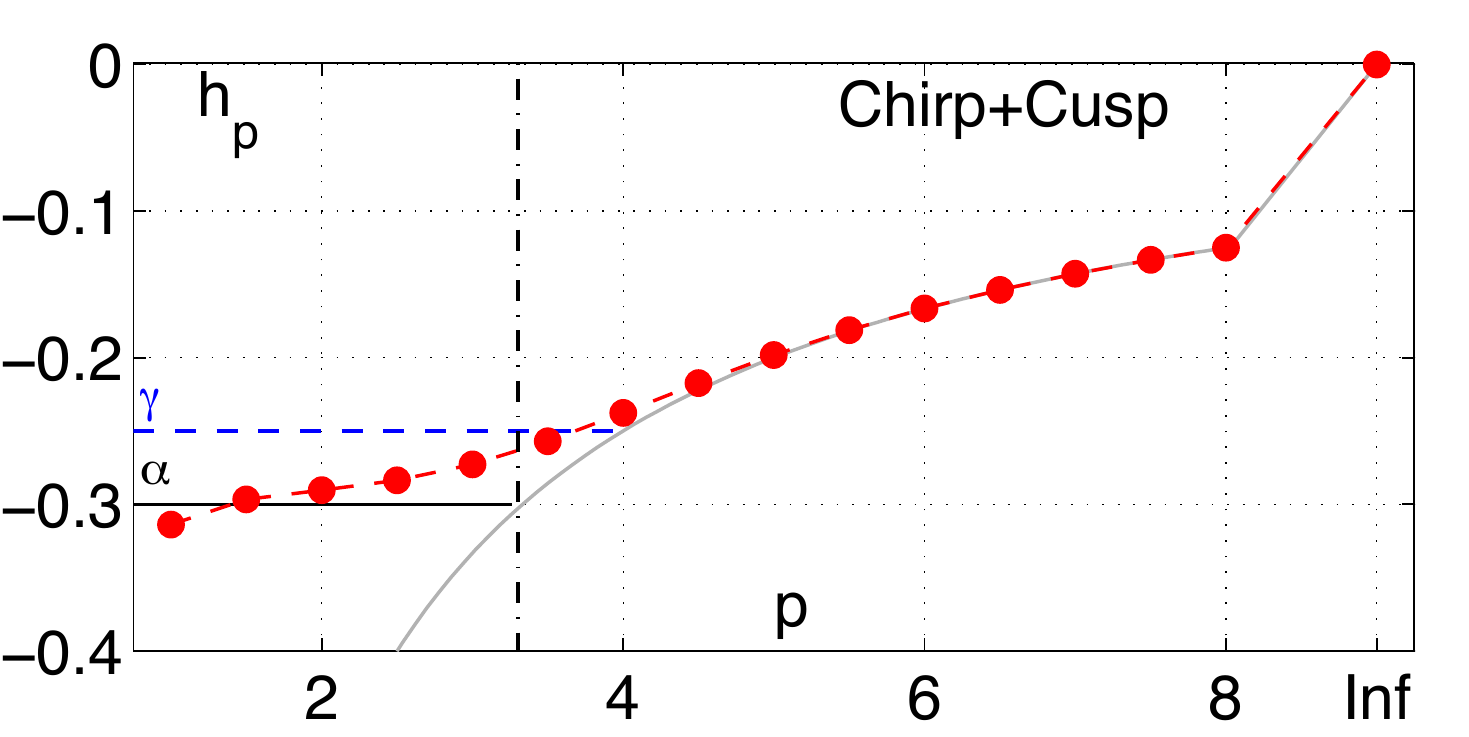}&
\includegraphics[height=30mm]{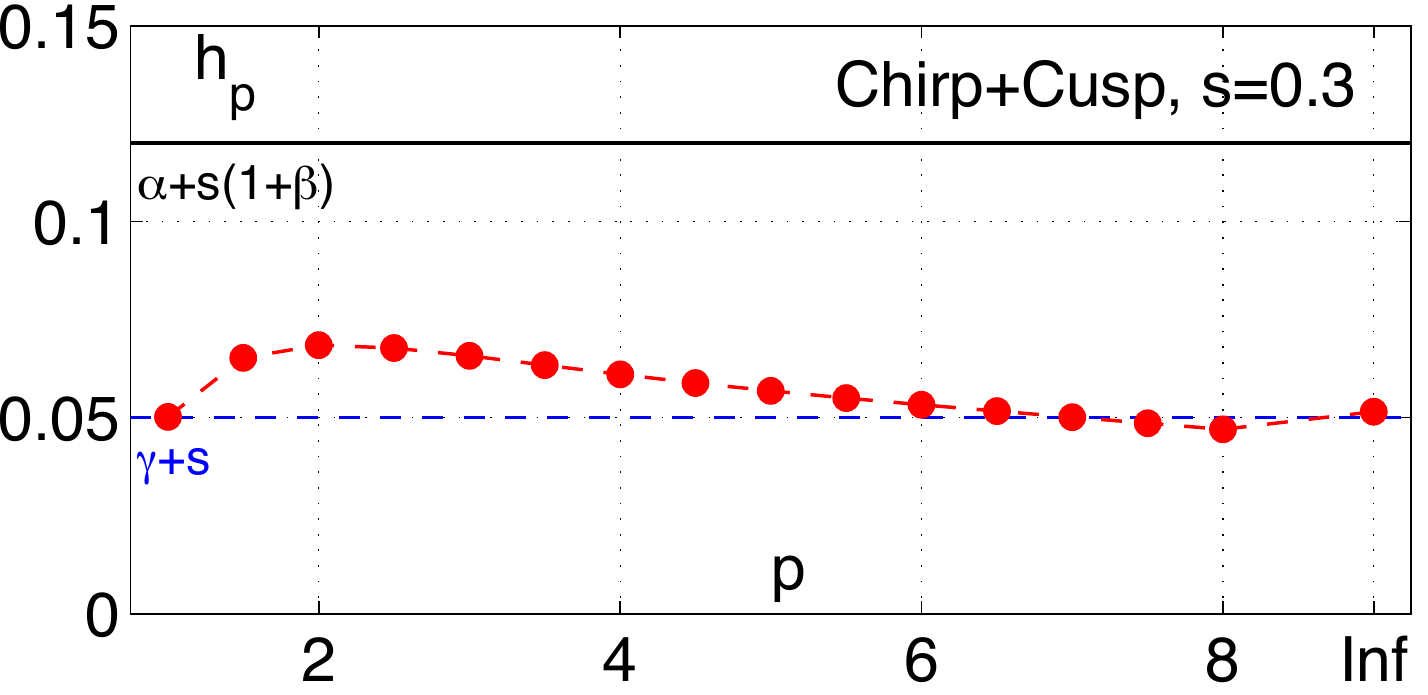}
\end{tabular}
\caption{\label{fig:cuspchirpFI} {\bf Chirp and cusp.} Sum of cusp and chirp singularities (top row). The $p$-leader based estimation of $p$-exponent (bottom row) for original functions (left column) and functions after a fractional integration of order $s=0.3$: theoretical value for $h_p(x_0)$ for cusp (blue dashed line) and chirp (black solid line), estimates (red disks) and $L^p$ limit (solid grey). The value for $p=\infty$ corresponds to the H\"older exponent. Only small values of $p$  allow to recover correctly the order of the dominant singularity: $\al = -0.3$}
\end{figure}

In contrast, we can use $p$-exponents with $p\leq -1/\alpha$ without a prior fractional integration.
The dominant singularity (chirp)  is $\alpha$ and $p$-exponents enable its correct estimation.
This example is illustrated numerically in Fig. \ref{fig:cuspchirpFI}. 
To conclude, in the presence of oscillating singularities, using  $p$-exponents should be preferred to the recourse to fractional integration, which induces smoothing and may make some of the singularities actually existing in data invisible.

\section{Conclusions and Perpectives}

As an alternative to the traditional H\"older exponent, the present contribution has put forward the use of a collection of new indices for pointwise, or local, regularity quantification: $p$-exponents, $p >0$, allowing notably the practical use of the notion of negative regularity, of significant importance for real-world data and applications.  
A corresponding wavelet framework enabling the practical estimation of $p$-exponents, the $p$-leaders, has been theoretically studied and numerically illustrated at work. 

The construction and analysis of $p$-exponents and $p$-leaders has yielded additional comprehensive understanding on local regularity analysis. 
First, it has shown that the choice of a particular index to quantify local regularities forces the selection of the relevant multiscale quantities for theoretical and practical analyses (while it is often misleadingly understood that the choices of the regularity index and of the multiscale quantities --- increments, wavelet coefficients --- can be independently achieved): 
While the H\"older exponent is intimately associated to oscillations and wavelet leaders,  $p$-exponents require the use of $p$-leaders. 
Second, the present work has shown that all regularity indices are not equivalent in general;
$p$-exponents hence permit to enrich the classification of singularities in terms of $p$-invariant versus non $p$-invariant, canonical versus oscillating, and oscillating balanced versus oscillating lacunary behaviors. 

The present contribution focused on the definitions and properties of $p$-exponents and $p$-leaders. 
It has hence been illustrated on isolated singularities for the sake of pedagogical exposition. 
The general aim of the present work however is multifractal analysis, aiming to study collections of intertwined singularities.
The construction of a $p$-leader based multifractal formalism is the topic of a companion article \cite{PART2}, which illustrates the additional benefits  of $p$-leaders even in situation where $p$-exponents and H\"older exponent coincide. 
This new formalism will also permit to show how $p$-leaders connect, generalize and outperform theoretically and practically \emph{Multifractal Detrended Fluctuation Analysis}, an earlier multifractal formalism of large popularity notably in biomedical applications (cf. \cite{kantelhardt2002multifractal} for the seminal definition). 

{\sc Matlab} routines designed by ourselves implementing $p$-leaders and the estimation procedures for $p$-exponents will be made publicly available to the research community at the time of publication of the present article.

\section*{Acknowledgements}

This work was supported by grants ANR AMATIS $\# 112432 $, 2012-2015, ANPCyT PICT-2012-2954 and PID-UNER-6136.

%\clearpage
%\newpage
\appendix

\section{\boldmath Proof of the conditions satisfied by the $p$-exponent function (first part of Theorem~\ref{th:tha})}
\label{sec:prooftha}

  { \bf Proof:} If $X \in L^p_{loc}$, then
  \[ \int_{B(0,a)} | X(u+x_0) |^p du \leq C \]
  so that the left hand-side of 
     (\ref{pexpa}) is bounded by $a^{-d/p}$, thus $X \in T^p_{-d/p} (x_0)$ and the first point of Theorem~\ref{th:tha} holds. 

  Assume that $X \in T^p_\alpha (x_0)$ for a $p >1$ and let $q$ be such that $1 < q  <p$.  We denote by $P$ the Taylor polynomial of $X$ for the $T^p_\alpha (x_0)$ regularity condition. Let $B= B(x_0, r)$,  $p'= p/q$ and $q'$ the conjugate exponent of $p'$, i.e.  $q'$ satisfies $\frac{1}{p'} + \frac{1}{q'}=1$. Using H\"older's inequality,  
   \[
\begin{array}{rl}
     \displaystyle\int_B | X (x) -P(x-x_0)|^q \; dx  =  &    \displaystyle\int_B | X (x) -P(x-x_0)|^q \; 1_B (x) \; dx  \\ & \\ 
 \leq  &    \left(   \displaystyle\int_B | X (x) -P(x-x_0)|^{qp'} \; dx \right)^{1/p'} \left(   \displaystyle\int    (1_B (x) )^{q'} \; dx \right)^{1/q'}   \\ & \\ 
  =  &       \left(   \displaystyle\int_B | X (x) -P(x-x_0)|^p \; dx \right)^{1/p'}  \left( C  r^d \right)^{1/q'}   \\  & \\ 
  \leq & C  \left(   r^{\al p +d} \right)^{1/p'}  r^{d/q'}  \\ & \\  =  & C r^{\al q +d}.  
\end{array}
\]
It follows that $X \in T^q_\alpha (x_0)$ and, up to the order $h_p (x_0)$,  the same Taylor polynomial can be used for $T^q_\alpha (x_0)$ as for $T^p_\alpha (x_0)$.  Hence the 
second point of Theorem~\ref{th:tha} holds, and also the assertion concerning the Taylor polynomial in Section \ref{sec:pexpodef}. \\ 

In order to prove the concavity of the function $ s \rightarrow h_{1/s}(x_0)$, we assume that  $X \in T^p_\alpha (x_0) \cap T^q_\beta (x_0)$. We pick $\gamma \in (0,1)$ and define $s$ by 
\[   \frac{1}{s} = \frac{\gamma}{p}  + \frac{1-\gamma}{q}  ; \]
 we have to prove that $X \in T^s_\delta (x_0)$ where 
\[ \delta = \gamma \al + (1-\gamma) \beta. \]
We pick for $P$ the Taylor polynomial of largest possible degree of $X$. Denote by $\parallel X \parallel_p$ the $L^p$ norm of $X$ on  $B(x_0, r)$; the assumptions imply that 
\[  \parallel X -P \parallel_p \leq C r^{\al + d/p}  \quad \mbox{ and} \quad \parallel X -P \parallel_q \leq C r^{\beta + d/q} . \]
  By interpolation, it follows that 
  \[  \parallel X -P \parallel_s \leq C r^{\gamma (\al + d/p) + (1-\gamma) (\beta + d/q)}  = C   r^{\gamma \al + (1-\gamma) \beta + d/s} , \]
  so that  $X \in T^s_\delta (x_0)$.  Hence the concavity of $ s \rightarrow h_{1/s}(x_0)$ follows.

\section{\boldmath Proof of the optimality of Theorem~\ref{th:tha}: Construction of functions with arbitrary $p$-exponents}
\label{sec:proofprop2}

We now show the optimality of Theorem ~\ref{th:tha}  by showing that the function  $ p \rightarrow h_{p}(x_0)$  can be any function satisfying the conclusions of Theorem ~\ref{th:tha}.    
The proof  consists of a constructive example.  \\

 For the sake of simplicity, we do the construction in dimension $d=1$;  the reader will easily check that it extends to arbitrary dimensions. 
    
    The functions that we will consider are parametrized by two  sequences $\theta (l)$ and $\omega (l)$ which satisfy the following conditions:
    
    \begin{enumerate}
    \item ${\omega (l) }+ {l }\rightarrow + \infty $
\item $\displaystyle\sum_{l=1}^\infty  2^{-\omega (l) -\theta (l)  }  < \infty$ 
\end{enumerate}

 The  function $F_{ \theta, \omega}$  is defined as follows (we take here $x_0 =0$) :
\[   
\begin{cases}
&\mbox{ if} \; x \in [2^{-l} , 2^{-l} +  2^{-\omega (l)}]   \mbox{ for an integer }  \; l >0,   \;  \mbox{ then} \qquad F_{ \theta, \omega} (x) = \frac{1}{l^2}2^{-\theta (l)} ,\\
&  \mbox{ otherwise} \qquad   F_{ \theta, \omega} (x) = 0. 
\end{cases}
\]  \vspace{3mm} 

Note that the first condition that we imposed implies a lacunarity on the construction: As $x \rightarrow 0$, $F_{ \theta, \omega}$ vanishes on a larger and larger proportion of points. 
The second  condition 
implies that $F_{ \theta, \omega}$  belongs to $L^1$ in the neighbourhood of 0. 

First, it is clear that $F_{ \theta, \omega}$ is  bounded in the neighbourhood of $x_0$ if and only if 
\BE \label{bounded}  \exists C \in \RR \quad \mbox{ such that} \qquad \forall l \geq 0, \qquad \theta (l) \geq C. \EE

\BL If (\ref{bounded}) holds, then the H\"older exponent of $F_{ \theta, \omega}$ at 0  is given by
\BE  \label{holexpthetome}  h_{\theta, \omega} (0) = \liminf_{l\rightarrow + \infty} \frac{\theta (l)}{l} . \EE 
\EL

{ \bf Proof:} Indeed, we  first note that $F_{ \theta, \omega}$ is continuous at $0$. Remark also that for $\alpha\geq 0$ the local maxima of the  function $F_{ \theta, \omega}(x) / x^\al  $  are obtained at the points $2^{-l}$. Thus if we pick $P =0$ in (\ref{equ:tpe}), then  computing the H\"older exponent amounts to find the supremum of the $\alpha$ such that one can find $C>0$ such that for all $l\geq 0 $ $\left|F_{ \theta, \omega}(2^{-l}) / 2^{-\al l}\right|\leq C $. This is easy to check that this supremum is exactly given by (\ref{holexpthetome}).\\

We still have to check that the estimate obtained by taking $P=0$ is the best possible. As regards the constant term, we separate two cases: if $\theta(l)$ does not tend to $+ \infty$ then   it follows  that $h_{\theta, \omega} (0) = 0$;  and if $\theta(l)\rightarrow + \infty$,  the constant term necessarily  vanishes.  For higher order terms, one argues by induction on the valuation of $P$, noticing that $F_{ \theta, \omega}$  vanishes  on the intervals  $[ 2^{-l} + 2^{-\omega ( l)} , 2^{-l+1} ] $,  which implies that, on such interval, the choice of  non-vanishing terms for  $P$ would lead to  worse estimates. \\ 

We now come back to the proof of the second part of Theorem ~\ref{th:tha} and  we no longer make the assumption that (\ref{bounded}) holds.   In order to estimate  the $p$-exponent of $F_{ \theta, \omega}$ at 0, we start by considering  the quantity 
\BE \label{holexpasetim}  \frac{1}{a} \int_{-a}^a | F_{ \theta, \omega}(x) |^p dx  \EE
(we thus take $P =0$ in the definition of the $T^p_\al $ regularity). 
First assume that $2^{-l_0} + 2^{-\omega (l_0) } \leq a \leq 2^{-l_0+1} $; 
then 
\BE\label{diver}  \frac{1}{a} \int_{-a}^a | F_{ \theta, \omega} (x) |^p dx  = \frac{1}{a}   \sum_{l = l_0}^\infty \frac{1}{l^{2p}}2^{-\omega (l) } 2^{-p\theta (l)  } . \EE
Therefore, the critical value $p(0)$ is given by 
\[ p(0) = \sup \{ p:  \mbox{  (\ref{diver}) converges} \} . \]

We will now make explicit choices for the sequences  $\theta (l)$ and $\omega (l)$. We pick the variable  $s = 1/p$, so that  the conditions satisfied by  the function \[ \rho (s) = \phi \left( \frac{1}{s} \right) \]   are that 
it  is a concave increasing function defined on $(1/p(0), 1]$ and  satisfying \[  \rho (s) \geq -\frac{1}{p(0)}.\] 
Concavity implies that  $\rho$ can  be obtained as the infimum of a countable family of affine functions 
\[ \rho_n (s) = a_n s + b_n \] 
which all satisfy
 \BE \label{rhon}   \forall s , \rho_n (s) \geq \rho (s). \EE
 We can also assume that there exists a dense sequence $s_n$ such that $\rho_n (s_n) = \rho (s_n)$.  
We now pick  the functions $\theta$ and $\omega$ as follows:
Note that any integer $l\geq 1$ can be written in a unique way under the form  $l = 2^n (2k +1)$, with $n \geq 0$ and $k \geq 0$;  then 
\[ \left\{ \begin{array}{rl}  \omega (l) =  & (a_n +1) l \\  \theta (l) =  & b_n l . 
 \end{array} \right. \] 
The result  follows  from the  fact that partial sums of the righthand side of (\ref{diver}) satisfy 
\BE \label{trbis}   \frac{1}{a}   \sum_{l = l_0}^m  \frac{1}{l^{2p}} 2^{-\omega (l) } 2^{-p\theta (l)  }  \leq   \frac{1}{a}   \sum_{l = l_0}^m  \frac{1}{l^{2p}} 2^{-( a_n +1+pb_n)l }   \EE
which, using  (\ref{rhon}) and $p \geq 1$, is bounded by 
\[  C 2^{l_0}   \sum_{l= l_0}^m  \frac{1}{l^2}  2^{-(p \rho (1/p)+1 )l } . \]
Using that the exponent is strictly positive (since $\phi$ is defined on $[1,p(0)] $), this sum is bounded by 
\BE \label{demorho}  C 2^{l_0}   2^{-(p \rho (1/p)+1 )l_0} \leq  C  \left( a^{\rho (1/p)}\right)^p ,\EE 
which is the required upper bound. 

If we now take $a$ in the interval $[2^{-l_0} , 2^{-l_0} + 2^{-\omega (l_0)}] $, the integral is bounded by the value at $2^{-l_0} + 2^{-\omega(l_0)}$ which has already been estimated, so that 
(\ref{tpx}) still holds.  We still have to check that the estimate obtained by taking $P=0$ is best possible. As in the H\"older case,  one  argues again  by induction on the valuation of $P$, noticing that  $F_{ \theta, \omega} $ vanishes  on the interval  $[ 2^{-l_0} + 2^{-\omega (l_0)}, 2^{-l_0} ] $, so that the order of magnitude of 
$ \int |F_{ \theta, \omega}  (x)  -P( x) |^p dx$  on such intervals will be given by the integral, on this interval,  of the first non-vanishing term of $P$. This remark immediately yields that, indeed, the choice $P= 0$ is  best possible.  \\

The lower bound is obtained by noticing that the first term of the sum in (\ref{trbis}) is larger than 
\[  \frac{1}{a}  \frac{1}{l_0^{2p}} 2^{-( a_n +1+pb_n)l_0 } \] 
which for $p = p_n =: 1/s_n$  and $a = 2^{-l_0} + 2^{-\omega (l_0)}$ is larger that 
\[  \frac{1}{2}  \frac{1}{l_0^{2p_n}} 2^{-p_n \rho (1/p_n)l_0 }; \] 
the lower bound is thus sharp at the points $p_n$.  Since they form a dense set, and since $\rho$ is continuous, (\ref{phip}) is proved. 

We have therefore obtained that the $p$-exponent of $F_{ \theta, \omega}  $ at the origin is   
\BE  \label{pexpfab} h_p (0) = \rho (1/p) ,  \EE  
and Theorem~\ref{th:tha}  is proved.

\section{Proof of Theorem~\ref{th:propo4bis} }
\label{sec:proofprop4} 

{ \bf Proof:}   We start by introducing a useful notation.  Let $p \in (0, + \infty)$; we will say that the sequence $\left( \pL_{\la_j (x_0)}\right)_{j\leq 0}$    satisfies  
\[ \left( \pL_{\la_j (x_0)}\right)_{j\leq 0}   \sim 2^{h j} \qquad \mbox{ at} \quad x_0 \]
 if the two following conditions hold:

\BE \label{upbouregpunp} \forall \ep >0, \quad \mbox{ for $j$ small enough,} \quad  \pL_{\la_j (x_0)} \leq 2^{(h - \ep)j}  \EE
\BE  \forall \ep >0, \quad \exists j_n \rightarrow -\infty \, : \qquad  \pL_{\la_{j_n} (x_0)} \geq 2^{(h + \ep)j_n} . \EE
Note that these two conditions can be rewritten as 
\BE \label{quatcingbis}  \liminf_{ j \rightarrow -\infty} \frac{\log \left(\pL_{\la_j (x_0)} \right) }{\log (2^{j}) } =h; \EE
and, since $\eta_X (p) >0$,  the wavelet characterization of $p$-exponents given by (\ref{carachqf}) exactly means that (\ref{quatcingbis}) is equivalent to $h_p (x_0) =h$. 
It follows from (\ref{carachqft}) that $X$ has  a canonical singularity  of exponent $h_{p_0}$ at $x_0$  if and only if 
\[ \mbox{ $\ell^{(p_0)}_{\la_j (x_0)}  \sim 2^{h_{p_0} j} $ at $x_0\quad $ and $\quad \exists s >0$: $\ell^{(p_0,s)}_{\la_j (x_0)} \sim 2^{(h_{p_0} +s )   j} $ at $x_0$. }\] 
Theorem~\ref{th:propo4bis}  will  be a direct consequence of   the following  result. 

\BP \label{prop:carcatpcusp} 
An admissible distribution $X$ has  a canonical singularity  of exponent $h_{p_0}$ at $x_0$   if and only if (\ref{upbouregpunp}) holds  for $p=p_0$ and if, for any $\ep >0$,  there exists a sequence of dyadic cubes $\la_n$ of scales $j_n \rightarrow -\infty$  such that 
\BE\label{eq:caractpcus} \left\{\begin{array}{ll} dist( x_0, \la_n  ) \leq 2^{(1-\ep) j_n } \\  \\ d_{\la_n}  \geq 2^{(h_{p_0}+\ep) j_n }  \end{array} \right.  
\EE 
where $d_{\lambda}$ is  the quantity defined in \eqref{pleaderssup}.
\EP

{ \bf Proof of Proposition \ref{prop:carcatpcusp}:} Let us assume that $X$ has  a canonical singularity of exponent $h_{p_0}$ at $x_0$.   We will first prove that

\
\BE \label{quat} \forall \ep >0, \quad \exists j_n \rightarrow -\infty : \quad 
 \left(  \sum_{ \lambda' \subset 3 \lambda_j (x_0) , \;  j_n\geq j' \geq (1+\ep) j_n }  d_{\lambda'} ^{p_0} 2^{d(j_n-j')} \right)^{1/{p_0}}  \geq 2^{(h_{p_0} + \ep ) j_n}   .\EE 

We prove it by contradiction. Indeed, if it were no true, then there would exist $\ep >0$ such that, 
\BE \label{eqcinq}  \forall j\leq 0 \qquad  \left(  \sum_{ \lambda' \subset 3 \lambda_j (x_0) , \;  j\geq j' \geq (1+\ep) j } d_{\lambda'}^{p_0} 2^{d(j-j')} \right)^{1/{p_0}}  \leq 2^{(h_{p_0} + \ep ) j} . \EE  

Let $s > 0$; consider the  corresponding $({p_0}, s)$-leader \eqref{ptleaders}:
\BE \label{pgammalead} \left(  \sum_{ \lambda' \subset 3 \lambda_j (x_0)  } (d_{\lambda'} 2^{s j'} )^{p_0} 2^{d(j-j')} \right)^{1/{p_0}} ; \EE
 we split the sum into two parts, depending whether $ j' \geq (1+\ep) j $ or $ j'  < (1+\ep) j $. 
 
It follows from (\ref{eqcinq}) that, in (\ref{pgammalead}), the term corresponding to  $j' \geq (1+\ep) j$ is bounded by $2^{(h_{p_0} + \ep + s  ) j} $. Consider now the
term corresponding to  $j'  <  (1+\ep) j$; since
$ 2^{ s   j'} \leq 2^{s (1+ \ep ) j},$
it follows that it is bounded by $2^{(h_{p_0} +  s (1+ \ep )  ) j} $. Therefore 
\[ h_{p_0, s} (x_0)  \geq 
h_{p_0} +  s + \ep (\min (1, s )  ) ,   \]
which contradicts the fact that $x_0$ is a canonical singularity. Therefore (\ref{quat}) holds. \\

We now prove that (\ref{quat}) implies that (\ref{eq:caractpcus}) holds.  Let $\ep'>0$. We first estimate the number of coefficients on which the sum in (\ref{quat}) bears using $\ep'$. 
At the first generation $j' = j_n-1$,  there are $2^d$ subcubes of $\la_{j_n} (x_0)$;   at each generation, each subcube of the  previous one is split into $2^d$ cubes; we go  down to the generation $j' = [ (1+\ep' )j_n ] =  j_n + [ \ep' j_n ] $, where there are $(2^d)^{[-\ep' j_n]} $ cubes. Therefore, the total number of cubes considered is bounded by $3 \cdot (2^d)^{[-\ep' j_n]} $. 

 It follows that  one of the terms of  
 the sum in (\ref{quat})  satisfies 
 \[  d_{\lambda'} ^{p_0} 2^{d(j_n-j') } \geq C 2^{(h_{p_0} + \ep' ) {p_0}j_n } 2^{d\ep' j_n }\] 
 which, using that $j_n \geq j' \geq j_n + \ep' j_n$,  implies that 
 \[ d_{\lambda'}  \geq C 2^{h_{p_0} j'  } 2^{C' \ep' j' } . \] 
 Furthermore, the conditions $\lambda' \subset 3 \lambda_{j_n} (x_0) $ and $j_n \geq j' \geq j_n+ \ep' j_n$ together  imply that the cube $\la'$ satisfies 
 the first condition of (\ref{eq:caractpcus}).\\ 
 
 We now prove the converse part in Proposition \ref{prop:carcatpcusp}. Assume that (\ref{upbouregpunp}) and (\ref{eq:caractpcus})   hold. We denote by $\tilde{\la}$ the smallest cube of the form $\la_j (x_0)$ such that $3 \tilde{\la}$ contains $\la'$. It follows from (\ref{eq:caractpcus}) that the scale $\tilde{j}$ of $\tilde{\la}$ satisfies $j' (1-\ep') \leq \tilde{j}\leq j'$. 
 Therefore the corresponding ${p_0}$-leader $\ell^{({p_0})}_{\tilde{\la}} $  satisfies 
 \BE \label{este1} L^{({p_0})}_{\tilde{\la}} \geq  \left( d_{\lambda'} ^{p_0} 2^{d(\tilde{j}-j') } \right)^{1/{p_0}}  \geq 2^{(h_{p_0} + \ep') j'} 2^{C \ep'\tilde{ j} } \geq 2^{h_{p_0} \tilde{j} } 2^{C \ep' \tilde{j}} ;\EE
  and, similarly, 
 the corresponding $({p_0}, s ) $-leader $\ell^{({p_0}, s )}_{\tilde{\la}} $  satisfies 
 \BE \label{este2} \ell^{({p_0},s)}_{\tilde{\la}} \geq  \left( (d_{\lambda'} 2^{s j' })^{p_0} 2^{d(\tilde{j}-j') } \right)^{1/{p_0}}  \geq 2^{(h_{p_0} + s + \ep') j'} 2^{C \ep' \tilde{j} } \geq 2^{(h_{p_0} + s) \tilde{j} } 2^{C \ep' \tilde{j}} . \EE
Remark that any $\ep>0$ can be written $C\ep'$, thus $h_{p_0, s} (x_0) \leq h_{p_0}+s$.  Note that (\ref{upbouregpunp}) implies that $h_{p_0} (x_0) \geq h_{p_0}$, and the lower bound $h_{p_0, s} \geq h_{p_0}+s$ then follows from general results on fractional integration, see \cite{JaffToul,JAFFARD:2010:A}. We have thus obtained that  
  $X$ has a canonical singularity   of exponent $h_{p_0}$ at $x_0$. \\ 
  
  We now turn to the proof of Theorem \ref{th:propo4bis}.  First, we note that the derivation of (\ref{este1})  also holds for $p < p_0$, and it follows that 
  \[ \forall p < p_0, \qquad h_p (x_0)  \leq h_{p_0}. \]
  Since  $ p\mapsto h_p (x_0)$   is decreasing,  it follows that $\forall p < p_0, $  $h_p (x_0) = h_{p_0}$.  Since (\ref{este1})  also holds for $p < p_0$, by  the same arguments as above, we deduce that 
  \[ \forall  p < p_0, \quad \forall s >0, \qquad h_{p, s} (x_0) = h_{p_0}+s. \]

\section{Proof of Proposition \ref{prop:ss}}
\label{proof:selfsimilar}
{ \bf Proof:}  We can assume that $x_0 =0$. We first prove the second part of the proposition; let us define the functions $ \omega_+   $ and  $ \omega_-  $  by (\ref{eq:prop3}); this implies that they belong to $L^p$ on $[0, \log a]$. The self-similarity condition exactly expresses the fact that they are $\log  a$ periodic functions. The condition $\al > -1/p$ then follows from the condition that $X \in L^p$ in the neighbourhood of 0.   

We now compute the $p$-exponent at $x_0$. We first assume that the Taylor polynomial  of $X$ vanishes at $0$.   For $r = 2^{-J}$, 
\[ \begin{array}{rl}
\displaystyle\frac{1}{r} \int_0^r |X|^p dx = & 2^J\displaystyle\int_0^{2^{-J}} | x|^{\al p} ( \ome_+ (\log x))^p dx  \\ & \\ 
= & 2^J \displaystyle\sum_{j= J+1}^\infty \int_{2^{-j}}^{2 \cdot 2^{-j}} | x|^{\al p} ( \ome_+ (\log x))^p dx   \\ 
& \\ 
= & 2^J\displaystyle\sum_{j= J+1}^\infty  2^{-j} \int_{1}^{2 } 2^{-\al p j} | u|^{\al p} ( \ome_+ (\log u))^p du  
 \end{array} \]
which is finite if and only if $\al p > -1$, in which case its value is 
$  C \cdot    2^{-\al p   J} $. The quantities are of the same order of magnitude if $r$ lies between two powers of type $2^{-J}$, and the computation is the same for $x <0$. Therefore the $p$-exponent of $X$  at $0$ is $\al$.

We now check that the Taylor polynomial  of $X$ vanishes at $0$.  First, if $ \al$ is not an integer, this is clear, since  a polynomial of degree less than $\al$ would bring a higher order contribution in the above computation. For the same reason, if $ \al$ is  an integer, then the polynomial necessarily is a monomial of order $\al$. Now, if $P (x) =  A x^\al$, then we note that $X-P$ also is self-similar of order $\al$, so that the above computation still applies, and the integral can be of a smaller order of magnitude only if  $X$ coincides with $P$ on both sides of the origin, which we excluded in the assumptions.

Let us now prove that $X$ has a canonical singularity at the origin. We note that,  through an integration by parts, one easily checks that,  if  $X$ is  a self-similar function of scaling ratio  $a$ and exponent $\al$ at $x_0$, then its primitive   which vanishes at 0 
is self-similar of scaling ratio  $a$ and exponent $\al +1$. It follows from the previous computations that its $p$-exponent is $\al +1$,  so that (\ref{equ-hpg}) holds with $s=1$, hence Proposition  \ref{prop:ss} is proved.

\bibliography{DFA2014}

\end{document}